\documentstyle[11pt,leqno,amsfonts,amssymb]{article}
\newtheorem{definition}{Definition}[section]
\newtheorem{lemma}[definition]{Lemma}
\newtheorem{theorem}[definition]{Theorem}
\newtheorem{proposition}[definition]{Proposition}
\newtheorem{corollary}[definition]{Corollary}
\newtheorem{remark}[definition]{Remark}

 1  

\newenvironment{proof}{\noindent{\bf Proof~:}}{\QED\medskip}
\def\QED{\hskip0.1em\hfill\null\ \null\nobreak\hfill
\kern3pt\lower1.8pt\vbox{\hrule\hbox
{\vrule\kern1pt\vbox{\kern1.7pt \hbox{$\scriptstyle
QED$}\kern0.2pt}\kern1pt\vrule}\hrule}}

\renewcommand{\a}{\alpha}

\newcommand{\ZZ}{{\mathbb{Z}}}

\newcommand{\RR}{{\mathbb{R}}}

\newcommand{\CP}{{\mathbb{CP}}}

\newcommand{\iso}{\cong}
\newcommand{\isom}{\stackrel{\cong}{\longrightarrow}}

\newcommand{\im}{\hbox{im}}
\renewcommand{\Re}{\hbox{Re}}
\renewcommand{\Im}{\hbox{Im}}
\newcommand{\rank}{\hbox{rank}\,}
\newcommand{\surj}{\twoheadrightarrow}
\newcommand{\inc}{\hookrightarrow}
\newcommand{\too}{\longrightarrow}

\newcommand{\PD}{\hbox{PD}}
\newcommand{\la}{\langle}
\newcommand{\ra}{\rangle}
\newcommand{\bd}{\partial}

\newcommand{\ox}{\otimes}
\newcommand{\bV}{\bigwedge V}
\newcommand{\bC}{\bigwedge C}
\newcommand{\hV}{\hat{V}}
\newcommand{\hC}{\hat{C}}
\newcommand{\hN}{\hat{N}}

\begin{document}

\title{\bf Formality of
       Donaldson submanifolds}

\author{Marisa Fern\'andez and Vicente Mu\~noz}

\date{October 29, 2003}

\maketitle
\begin{abstract}
We introduce the concept of $s$--{\em formal minimal model} as an
extension of formality. We prove that any orientable compact
manifold $M$, of dimension $2n$ or $(2n-1)$, is formal if and only
if $M$ is $(n-1)$--formal. The formality and the hard Lefschetz
property are studied for the Donaldson submanifolds of  symplectic
manifolds constructed in~\cite{Do}. This study permits us to show
an example of a Donaldson symplectic submanifold of dimension
eight which is formal simply connected and does not satisfy the
hard Lefschetz theorem.
\end{abstract}

\section{Introduction} \label{introduction}
A symplectic manifold $(M,\omega)$ is a pair consisting of a
$2n$--dimensional differentiable manifold $M$ with a closed
$2$--form $\omega$ which is non-degenerate (that is, $\omega^n$
never vanishes). The form $\omega$ is called symplectic. By the
Darboux theorem, in canonical coordinates, $\omega$ can be
expressed as $\omega = \sum\limits_{i=1}^{n} dx^i\wedge dx^{n+i}$.
Therefore any symplectic manifold admits an almost complex
structure $J$ compatible with the symplectic form, which means
that $\omega(X,Y)= \omega(JX,JY)$ for any $X$, $Y$ vector fields
on $M$.

The simplest examples of symplectic manifolds are K\"ahler
manifolds; for example, the complex projective space $\CP^n$ with
the standard K\"ahler form $\omega_0$ defined by its natural
complex structure and the Fubini--Study metric. Gromov~\cite{Gr1,
Gr2} and Tischler~\cite{Ti} prove that if $M$ is a compact
symplectic manifold, of dimension $2n$, with an integral
symplectic form $\omega$, then there is a symplectic embedding $f
\colon (M,\omega)  \longrightarrow (\CP^{2n+1},\omega_0)$ such
that $f^*(\omega_0)=\omega$.

If $(M,\omega)$ is a compact symplectic manifold, then the de Rham
cohomology groups $H^{2k}(M)$, $k\leq n$, are non-trivial. (We
denote by $H^*(M)$ the real or de Rham cohomology, and by
$H^*(M;\ZZ)$ the integral one.) The problem of how compact
symplectic manifolds differ topologically from K\"ahler manifolds
led during the last years to the introduction of several geometric
methods for constructing symplectic manifolds. They include
compact nilmanifolds~\cite{BG, topo, Th}, symplectic blow
ups~\cite{McD1}, and fiber connected sums~\cite{Go}. The
symplectic manifolds there presented do not admit a K\"ahler
metric since either they are not formal or do not satisfy hard
Lefschetz theorem, or they fail both properties of compact
K\"ahler manifolds.

Let $(M,\omega)$ be a compact symplectic manifold of dimension
$2n$ with $[\omega]\in H^2(M)$ having a lift to an integral
cohomology class $h$. In \cite{Do} Donaldson proves the existence
of some integer number $k_0$ such that for any $k\geq k_0$ there
exists a symplectic submanifold $Z\inc M$ of dimension $2n-2$ that
realizes the Poincar\'e dual of $kh$, that is, $\PD [Z]=kh\in
H^2(M;\ZZ)$. We shall call these manifolds {\em Donaldson
symplectic submanifolds\/} (or, indistinctly, {\em Donaldson
submanifolds\/}) of $M$. Such manifolds satisfy a {\em Lefschetz
theorem on hyperplane sections\/}. This means that the inclusion
$j\colon Z\inc M$ is $(n-1)$--connected, i.e., up to homotopy $M$
is constructed out of $Z$ by attaching cells of dimension $n$ and
higher. In particular,
\begin{itemize}
 \item for $i<n-1$ there is an isomorphism $j^*\colon  H^i(M)\to
  H^i(Z)$ induced by $j$ on cohomology;
 \item for $i=n-1$ there is a monomorphism $j^*\colon H^i(M)\inc H^i(Z)$.
\end{itemize}

Our purpose in this note is to study the formality and the hard
Lefschetz theorem for Donaldson submanifolds of  symplectic
manifolds. As a consequence of this study, we prove that a compact
simply connected symplectic $8$-manifold constructed in
\cite{IRTU} is formal, while it does not satisfy the hard
Lefschetz theorem.

The description of a minimal model for a Donaldson submanifold of
a symplectic manifold can be very complicated even for the degree
$n-1$. This is the reason for which we need first to weaken the
condition of formal manifold to {\em s--formal manifold} $(s\geq
0)$ as follows. Let $(\bigwedge V,d)$ be a minimal model of a
differentiable manifold $M$ (of arbitrary dimension). We say that
$(\bigwedge V,d)$ is a {\em s--formal minimal model}, or $M$ is a
{\em s--formal manifold}, if for each $i\leq s$ the subspace $V^i$
of $V$, consisting of the generators of degree $i$, decomposes as
a direct sum $V^i= C^i\oplus N^i$ where the spaces $C^i$ and $N^i$
satisfy the three following conditions:

\begin{enumerate}
\item $d(C^i) = 0$,
\item the differential map $d\colon
N^i\longrightarrow \bigwedge V$ is injective,
\item any closed element in the ideal
$I(\bigoplus\limits_{i\leq s} {N^i}) = N^{\leq s}\cdot
(\bigwedge V^{\leq s}))$, generated by $\bigoplus\limits_{i\leq
s} {N^i}$ in $\bigwedge (\bigoplus\limits_{i\leq s} {V^i})$, is
exact in $\bigwedge V$.
\end{enumerate}

Note that if $M$ is $(s+1)$--formal, then $M$ is $s$--formal. All
connected manifolds are obviously $0$-formal, and there are
examples of non-formal manifolds. For any $s\geq 0$, we show
examples of manifolds which are $s$--formal but not
$(s+1)$--formal (see Example $5$ in Section~\ref{examples}). The relation
between
the $s$--formality and the formality is given in
Theorem~\ref{criterio2} of Section~\ref{maintheorem}. There we
prove that any orientable compact connected manifold $M$, of
dimension $2n$ or $(2n-1)$, is formal if and only if $M$ is
($n-1$)--formal. This means that the formality of $M$ is contained
in the ($n-1$) first subspaces $V^i$ ($1 \leq i \leq {n-1}$) of
the minimal model of $M$, and so we can ignore the other
subspaces.

As a consequence of Theorem~\ref{criterio2} we get Miller's
theorem~\cite{Mi} for the formality of a $(k-1)$--connected
compact manifold of dimension less than or equal to $(4k-2)$.
We show that any simply connected oriented compact
manifold $M$ of arbitrary dimension is $2$--formal, as well as if
$M$ has dimension $7$ or $8$, $M$ is formal if and only if is
$3$--formal.

In Theorem~\ref{equivalence2} we prove that if $M$ is a
$(n-2)$--formal compact symplectic manifold of dimension $2n$, and
$Z\inc M$ is a Donaldson submanifold, then $Z$ is formal.
Therefore, it can happen that $Z$ is formal but $M$ is non-formal.
Moreover, in Theorem~\ref{Lefschetz-Donaldson} we get the conditions
on $M$ under which it is possible to state that $Z$ satisfies the
hard Lefschetz theorem.

The paper is structured as follows. In Section~\ref{seccionformal}
we establish the concept of $s$--formal manifold. For such a
manifold $M$ we show, in Lemma~\ref{topologia}, that if $(\bV,d)$ is
the minimal model of $M$ then the minimal model $(\bigwedge W,d)$ of
the differential algebra $(H^*(M),d=0)$ is given by $(\bigwedge
V^{\leq s},d)$ by adding spaces $W^{> s}$ with suitable
differentials. The relation between $s$--formality and Massey
products is proved in Lemma~\ref{Massey}.

In Section~\ref{maintheorem} we determine the smallest value of
$s$ for which the $s$--formality is equivalent to the formality of
$M$ proving Theorem~\ref{criterio2}. Some consequences are
discussed; in particular, Miller's theorem~\cite{Mi}.
Section~\ref{sectionlefschetz} is dedicated to compact symplectic
manifolds $(M,\omega)$ with the $s$--Lefschetz property ($s\leq
{(n-1)}$), i.e., satisfying that the cup product
 $$
 [\omega]^{n-i}: H^i(M) \too H^{2n-i}(M)
 $$
is an isomorphism for all $i\leq s$. In
Section~\ref{sectionsubvarities} we prove that a Donaldson
submanifold $Z\subset M$ is hard Lefschetz if and only if $M$ has
the $(n-2)$--Lefschetz property. We also show the
Theorem~\ref{equivalence2} previously mentioned on the formality
of Donaldson submanifolds.

Finally, Section~\ref{examples} is devoted to the discussion of
examples illustrating the concepts and results of the previous
sections. Such examples reveal the existence of Donaldson
symplectic submanifolds satisfying one of the following
properties: formal and hard Lefschetz; neither formal nor hard
Lefschetz; or formal but not hard Lefschetz. Furthermore, some of
those are $4$--dimensional symplectic manifolds that have neither
K\"{a}hler metrics nor complex structures.

\section{$s$--formality and real homotopy}\label{seccionformal}

In this section, we establish the concept of $s$--formal minimal
model and show some properties for such manifolds. First, we need
some definitions and results about minimal models.

Let $(A,d)$ be a {\it differential graded algebra} (in the sequel,
we shall say just a differential algebra), that is, $A$ is a
graded commutative algebra over a field $K$, of characteristic
zero, with a differential $d$ which is a derivation, i.e.
$d(a\cdot b) = (da)\cdot b +(-1)^{\deg (a)} a\cdot (db)$, where
$\deg(a)$ is the degree of $a$. Throughout this article all vector
spaces and algebras are defined over the field $\RR$ of real
numbers unless there is indication to the contrary.

A differential algebra $(A,d)$ is said to be {\it minimal\/} if:
\begin{enumerate}
 \item $A$ is free as an algebra, that is, $A$ is the free
 algebra $\bigwedge V$ over a graded vector space $V=\oplus V^i$, and
 \item there exists a collection of generators $\{ a_\tau,
 \tau\in I\}$, for some well ordered index set $I$, such that
 $\deg(a_\mu)\leq \deg(a_\tau)$ if $\mu < \tau$ and each $d
 a_\tau$ is expressed in terms of preceding $a_\mu$ ($\mu<\tau$).
 This implies that $da_\tau$ does not have a linear part, i.e., it
 lives in $\bV^{>0} \cdot \bV^{>0} \subset \bV$.
\end{enumerate}

Morphisms between differential algebras are required to be degree
preserving algebra maps which commute with the differentials.
Given a differential algebra $(A,d)$, we denote by $H^*(A)$ its
cohomology.
$A$ is {\em connected\/} if $H^0(A)=\RR$, and $A$ is {\em
one--connected\/} if, in addition, $H^1(A)=0$.

We shall say that $({\cal M},d)$ is a {\it minimal model} of the
differential algebra $(A,d)$ if $({\cal M},d)$ is minimal and there
exists a morphism of differential graded algebras $\rho\colon
{({\cal M},d)}\longrightarrow {(A,d)}$ inducing an isomorphism
$\rho^*\colon H^*({\cal M})\longrightarrow H^*(A)$ on cohomology.

In~\cite{H} Halperin proved that any connected differential algebra
$(A,d)$ has a minimal model unique up to isomorphism. For $1$--connected
differential algebras, a similar result was proved by Deligne, Griffiths,
Morgan and Sullivan~\cite{DGMS,GM,Su}.


A minimal model $({\cal M},d)$ is said to be {\it formal} if there is a
morphism of differential algebras $\psi\colon {({\cal M},d)}\longrightarrow
(H^*({\cal M}),d=0)$ that induces the identity on cohomology. The
formality of a minimal model can be characterized as follows.

\begin{theorem}~\cite{DGMS}.\label{criterio1}
A minimal model $({\cal M},d)$ is formal if and only if we can write
${\cal M}=\bV$ and the space $V$ decomposes as a direct sum $V= C\oplus
N$ with $d(C) = 0$, $d$ is injective on $N$ and such that every
closed element in the ideal $I(N)$ generated by $N$ in $\bigwedge
V$ is exact.
\end{theorem}

A {\it minimal model\/} of a connected differentiable manifold $M$ is a
minimal model $(\bigwedge V,d)$ for the de Rham complex $(\Omega
M,d)$ of differential forms on $M$. If $M$ is a simply connected
manifold, the dual of the real homotopy vector space
$\pi_i(M)\otimes \RR$ is isomorphic to $V^i$ for any $i$. This
relation also happens when $i>1$ and $M$ is nilpotent, that is,
the fundamental group $\pi_1(M)$ is nilpotent and its action on
$\pi_j(M)$ is nilpotent for $j>1$ (see~\cite{DGMS,GM}).

We shall say that $M$ is {\it formal\/} if its minimal model is
formal or, equivalently, the differential algebras $(\Omega M,d)$
and $(H^*(M),d=0)$ have the same minimal model. (For details
see~\cite{DGMS,GM,L} for example.) Therefore, if $M$ is formal and
simply connected, then the real homotopy groups $\pi_i(M)\otimes
\RR$ are obtained from the minimal model of $(H^*(M),d=0)$.

From now on, we will consider only 
connected differentiable manifolds. 
In order to obtain some information on the formality of a
manifold, we introduce the concept of $s$--formality as follows.

\begin{definition}\label{primera}
Let $({\cal M},d)$ be a minimal model. We say that
$({\cal M},d)$ is $s$--formal
($s\geq 0$) if we can write ${\cal M}=\bV$ such that for each $i\leq s$
the space $V^i$ of generators of degree $i$ decomposes as a direct
sum $V^i=C^i\oplus N^i$, where the spaces $C^i$ and $N^i$ satisfy
the three following conditions:
\begin{enumerate}
\item $d(C^i) = 0$,
\item the differential map $d\colon N^i\longrightarrow \bigwedge V$ is
injective,
\item any closed element in the ideal
$I_s=I(\bigoplus\limits_{i\leq s} N^i)$, generated by the space
$\bigoplus\limits_{i\leq s} N^i$ in the free algebra $\bigwedge
(\bigoplus\limits_{i\leq s} V^i)$, is exact in $\bigwedge V$.
\end{enumerate}
\end{definition}

In what follows, we shall write $N^{\leq s}$ and $\bigwedge
V^{\leq s}$ instead of $\bigoplus\limits_{i\leq s}N^i$ and
$\bigwedge(\bigoplus\limits_{i\leq s}V^i)$, respectively. In
particular, $I_s=N^{\leq s}\cdot (\bigwedge V^{\leq s})$.

\begin{remark}
We must note that the conditions in Definition \ref{primera} are
not the same as to ask that $(\bigwedge V^{\leq s},d)$ is formal,
since in the later case, one should ask that every closed element
in the ideal $I(N^{\leq s})$ is exact in $\bigwedge V^{\leq s}$.
Moreover, Definition \ref{primera} implies that if ($\bigwedge
V$,d) is formal then it also is $s$--formal for any $s\geq 0$.
\end{remark}

\begin{definition}\label{primeravar}
A differentiable manifold $M$ is $s$--formal if its minimal model
is $s$--formal (in the sense of the previous Definition).
\end{definition}

\begin{remark}\label{rem:cw}
  The concept of $s$--formality can be defined for CW--complexes
  which have a minimal model $(\bV,d)$. Such a minimal model is
  constructed as the minimal model associated to the differential
  complex of piecewise-linear rational polynomial forms \cite{GM}.
\end{remark}

Let $M$ be an $s$--formal manifold with minimal model $(\bigwedge
V,d)$ as in Definition \ref{primera}. Clearly, the space
$\bV^{\leq s}$ has the decomposition
 \begin{equation} \label{eqn:v0}
  \bV^{\leq s} = (\bigwedge C^{\leq s}) \oplus
  N^{\leq s}\cdot  (\bigwedge V^{\leq s}).
 \end{equation}
For degree $i\leq s$ it is clear that $(\bigwedge V^{\leq
s})^i = (\bigwedge V)^i$. Then for $i\leq s$ we have a surjection
 $$
 (\bigwedge C^{\leq s})^i \surj H^i(M).
 $$

Before going into the study of $s$--formal manifolds, we show
examples of compact connected manifolds which are
$0$--formal but not $1$--formal. The simplest examples are the compact
nilmanifolds non-tori. Let $G$ be a connected rational nilpotent
Lie group, and denote by $\Gamma$ a discrete subgroup of $G$ such
that the quotient space $M=\Gamma{\backslash}G$ is
compact~\cite{Ma}. Such a manifold $M$ is called a compact
nilmanifold. Hasegawa's theorem~\cite{Ha} states that the tori are
the only formal compact nilmanifolds. We reformulate that theorem
as follows.

\begin{lemma}~\cite{Ha}.\label{formalitynilm}
Let $M =\Gamma{\backslash}G$ be a compact nilmanifold. Obviously
$M$ is $0$--formal. Moreover, $M$ is $1$--formal if and only if
$M$ is a torus, and so formal.
\end{lemma}

\begin{proof}
It is clear that any differentiable manifold is $0$--formal.
Consider $M =\Gamma{\backslash}G$ a compact nilmanifold. Then, a
minimal model $({\cal M},d)$ of $M$ is given by $(\bigwedge
({\frak g}^*),d)$, where ${\frak g}$ is the Lie algebra of $G$,
and $d$ is the Chevalley--Eilenberg differential in $\bigwedge
({\frak g}^*)$. Since all the generators of $\cal M$ have degree
$1$, according to Definition~\ref{primera} and
Definition~\ref{primeravar}, we get that $\cal M$ is $1$--formal
if and only if $\cal M$ is formal, and so $M$ is formal. This
completes the proof using Hasegawa's theorem.
\end{proof}

In
Section~\ref{examples} we construct examples of
{\it non-symplectic} manifolds which are $s$--formal but not
($s+1$)--formal for $s\geq 2$ (see Example $5$), and
we show examples of compact {\it symplectic} manifolds
which are $s$--formal but not
($s+1$)--formal for $s=1$ (see Example $3$) and for $s=3$ (see
Example $4$).

\bigskip
Next, we show the first properties of $s$--formal manifolds. For
such a manifold, the analogous result to Theorem~\ref{criterio1}
is the following lemma.


\begin{lemma}\label{topologia}
Let $M$ be a manifold with minimal model $(\bV,d)$. Then $M$ is
$s$--formal if and only if there is a map of differential algebras
 $$
 \varphi: (\bV^{\leq s},d) \too (H^*(M),d=0),
 $$
such that the map $\varphi^*: H^*(\bV^{\leq s},d) \too H^*(M)$
induced on cohomology is equal to the map $\imath^*:H^*(\bV^{\leq
s},d) \too H^*(\bV,d)=H^*(M)$ induced by the inclusion $\imath:
(\bV^{\leq s},d) \too (\bV,d)$.

In particular, $\varphi^*: H^i(\bV^{\leq s}) \too H^i(M)$ is an
isomorphism for $i\leq s$, and a monomorphism for $i=s+1$. So, if
$M$ is simply connected, then the dual of the real homotopy vector
space $\pi_i(M)\otimes \RR$ is isomorphic to $V^i=W^i$ for any
$i\leq s$, $(\bigwedge W,d)$ being the minimal model of
$(H^*(M),d=0)$.
\end{lemma}

\begin{proof}
Since $(\bV,d)$ is a minimal model of $M$, we know that there is a
morphism $\rho\colon {(\bV,d)}\longrightarrow {(\Omega M,d)}$
inducing an isomorphism $\rho^*$ on cohomology. Thus to prove the
{\em only if\/} part, it is sufficient to show that there is a map
of differential algebras
 $$
 \psi: (\bV^{\leq s},d) \too (H^*(\bV),d=0),
 $$
such that the map $\psi^*: H^*(\bV^{\leq s}) \too H^*(\bV)$
coincides with the map $\imath^*:H^*(\bV^{\leq s}) \too H^*(\bV)$
induced by the inclusion $\imath: (\bV^{\leq s},d) \too (\bV,d)$.
Then the map $\varphi$ given by $\varphi=\rho^*\circ \psi$ satisfies
the conditions that we need.

We define $\psi(x)=[x]$ for $x\in C^i$ and $\psi(x)=0$ for $x\in
N^i$, $i\leq s$. We extend $\psi$ to an algebra map $\psi: \bV^{\leq
s} \too H^*(\bV)$ by multiplicativity. We see that $s$--formality
implies that $\psi$ commutes with the differentials, as follows. Let
$x \in \bV^{\leq s}$. Then, from~(\ref{eqn:v0}), it follows that
$dx$ decomposes $dx=a+b$ with $a\in \bigwedge C^{\leq s}$ and $b \in
N^{\leq s}\cdot  (\bigwedge V^{\leq s})$. Since $a$ is closed, so is
$b=dx-a$ and hence exact by $s$--formality. Therefore $a=dx-b$ is
exact as well, and $\psi(dx)=\psi(a)=[a]=0$. This shows that $\psi:
(\bV^{\leq s},d) \too (H^*(\bV),d=0)$ is a map of differential
algebras.

Moreover, if $x\in\bV^{\leq s}$ is closed, decompose $x=a+b$ with
$a\in \bigwedge C^{\leq s}$, $b \in N^{\leq s}\cdot (\bigwedge
V^{\leq s})$. Then $b=x-a$ is closed, hence exact by $s$--formality.
So $\psi^*[x]=[\psi(x)]=\psi(x)=\psi(a)=[a]=[a+b]=[x]$, hence
$\psi^*=\imath^*$, as required.

To prove the {\em if\/} part, suppose that we have a map $\varphi:
(\bV^{\leq s},d) \too (H^*(M),d=0)$ satisfying $\varphi^*=\imath^*$.
We want to find an $s$--formal model for $M$, i.e.,\ $\bV=\bigwedge
\hV$ such that $\hV^i=\hC^i \oplus \hN^i$ satisfies the conditions
of Definition~\ref{primera}. Moreover we construct this model in
such a way that $\bV^{\leq i}=\bigwedge \hV^{\leq i}$ for all $i$.
Let us do this by induction on $i$. Suppose that we have defined
$\hV^{<i}=\hC^{<i} \oplus \hN^{<i}$, with $i\leq s$. Then we define
   $$
    N^i=\ker (\varphi: V^i \to H^i(M)/ \im (\varphi\colon
    \bigwedge V^{<i} \to H^i(M))).
   $$
For $x \in N^i$, let $a_x\in\bV^{<i}$ be a closed element such
that $\varphi(x)=[a_x]$ and set $\hat{x}= x-a_x$. This defines a
space $\hat{N}^i$ isomorphic to $N^i$. Consider the space $C^i$ given by
   $$
   C^i=\ker (d: V^i \to \bV/d (\bV^{<i})).
   $$
Now for $y\in C^i$, let $b_y\in \bV^{<i}$ such that $d y=db_y$.
This $b_y$ is well-defined up to a closed element, so we may
suppose that $\varphi(b_y)=0$. Set $\hat{y}=y-b_y$. This defines a
space $\hat{C}^i$ isomorphic to $C^i$. Now $\varphi(\hat{N}^i)=0$
and $d(\hat{C}^i)=0$. If we check that $V^i=C^i\oplus N^i$ then it
follows that $\hV^i=\hC^i\oplus \hN^i$ is isomorphic to $V^i$.

First, if $x\in N^i\cap C^i$ then $\varphi(x-a_x)=0$ and
$d(x-b_x)=0$. So $x-a_x-b_x$ is closed and $\varphi(x-a_x-b_x)=0$.
Therefore $x-a_x-b_x$ is exact, which contradicts the minimality
of the model. Thus $N^i\cap C^i=0$. Second, if $x\in V^i$ then
consider $\varphi(x)=[t]$ where $t\in V^i\oplus \bV^{<i}$ and $d
t=0$. Decompose $t=t_1+t_2$ where $t_1\in V^i$, $t_2\in \bV^{<i}$,
so that $t_1\in C^i$. Now $\varphi(x-t)=0$ so that $x-t_1\in N^i$.
Therefore $x\in C^i\oplus N^i$. The properties of Definition
\ref{primera} are now easy to verify for $\hV^i$.

For the final assertion, let $(\bigwedge W,d)$ be the minimal model
of $H^*(M)$. This can be constructed starting with $(\bV^{\leq
s},d)$, and with the map $\varphi: (\bV^{\leq s},d) \too
(H^*(M),d=0)$, and adding subspaces $W^{>s}$ with suitable
differentials (see~\cite{TO}).
\end{proof}

\begin{remark}
  The concept of $1$--formality appears in~\cite[Chapter 2]{Am}
  defined by the formulation given in Lemma~\ref{topologia}.
  In~\cite{Am} the $1$--formality is studied in connection with the
  fundamental group of the manifold.
\end{remark}

It is well known that all Massey products vanish for any formal manifold.
The relation between $s$--formality and Massey products is given
in the following lemma.

\begin{lemma}\label{Massey}
Let $M$ be an  $s$--formal manifold. Suppose that there are
cohomology classes $\alpha_i\in H^{p_i}(M)$, $p_i>0$, $1\leq i\leq
t$, such that the Massey product $\la
\alpha_1,\alpha_2,\dots,\alpha_t\ra$ is defined. If $p_1+p_2+
\cdots+p_{t-1}\leq {s+t-2}$ and $p_{2}+ \cdots+p_t\leq {s+t-2}$,
then $\la \alpha_1,\alpha_2,\dots,\alpha_t\ra$ vanishes.
\end{lemma}

\begin{proof}
Let $(\bV,d)$ be a minimal model of $M$. There exists a morphism
$\rho\colon {(\bV,d)}\longrightarrow {(\Omega M,d)}$ inducing an
isomorphism $\rho^*$ on cohomology. For $1\leq i\leq t$, denote by
$[a_i]\in H^{p_i}(\bV)$ the cohomology classes such that
$\rho^*[a_i]=\alpha_i$. To prove the Lemma we see that the Massey
product $\la [a_1],[a_2],\dots,[a_t]\ra$ vanishes.

First we show it for triple Massey products. Consider a Massey
product $\la [a_1],[a_2],[a_3]\ra$, $[a_i]\in H^{p_i}(\bV)$, with
$p_1+p_2\leq s+1$ and $p_2+p_3\leq s+1$. Suppose that $a_1\cdot
a_2=d \xi_1$ and $a_2\cdot a_3=d \xi_2$. Since the degree of
$\xi_j$ does not exceed $s$, we have that $\xi_j \in \bV^{\leq
s}$. Projecting onto the second summand of (\ref{eqn:v0}) we can
suppose that $\xi_j\in I_s$. By the $s$--formality,
  $$
  a_1 \cdot \xi_2+(-1)^{p_1+1} \xi_1 \cdot a_3
  $$
is exact, and so the triple Massey product $\la
[a_1],[a_2],[a_3]\ra$ vanishes.

The case of the higher Massey product is similar. Let us first
recall the definition (see~\cite{K,Mas,Ta}). If the Massey product
$\la [a_1],[a_2],\dots,[a_t]\ra$ is defined, then there are
elements $a_{i,j}$ of the minimal model $\bV$ of $M$, with $1\leq
i\leq j\leq t$, except for the case $(i,j)=(1,t)$, such that
$a_{i,i}$ is a cocycle representing $[a_i]$ and $d\,a_{i,j}=
\sum\limits_{k=i}^{j-1} {\bar a}_{i,k}\cdot a_{k+1,j}$, where
$\bar a=(-1)^{\deg(a)} a$. Then the Massey product $\la
[a_1],[a_2],\dots,[a_t]\ra$ is the set of all possible cohomology
classes of degree $p_1+ \cdots +p_t -(t-2)$ whose representatives
are $\sum\limits_{k=1}^{t-1} {\bar a}_{1,k} \cdot a_{k+1,t}$. If
one of these representatives is exact, then the Massey product
$\la [a_1],[a_2],\dots,[a_t]\ra$ is zero.

Now, for $k=1$, $a_{1,1}$ is a closed element representing the
cohomology class $[a_1]$ and we know that $\deg(a_1) = p_1 \leq
s$. For $k=2$, $d\,a_{1,2}=a_1 \cdot a_2$, that is, $\deg(a_{1,2})
= p_1+p_2-1$ which is  $\leq s$. 
For any $3\leq k\leq {t-1}$, $d\, a_{1,k}$ is a representative of
the Massey product $\la [a_1],[a_2],\dots,[a_k]\ra$. Thus,
$\deg(a_{1,k})=p_1+\cdots + p_k - (k-1)$ which is less than or
equal to $s$ by hypothesis. Hence $a_{1,1}\in \bC^{\leq s}$ and
$a_{1,k}\in \bV^{\leq s}$ for $2 \leq k\leq {t-1}$. In a similar
way we see that $\deg(a_{k+1,t})\leq s$ for $1 \leq k\leq {t-2}$;
and, for $k=t-1$, the element $a_{k+1,t}=a_{t,t}$ is a
representative of $[a_t]$ and so it has degree $\leq s$.
Therefore, $a_{1,t-1}$, $a_{2,t}$, $a_{1,k}$ and $a_{k+1,t}\in
\bV^{\leq s}$ for $2\leq k\leq {t-2}$, and $a_{1,1}$, $a_{t,t}\in
\bC^{\leq s}$. Using (\ref{eqn:v0}) one can project onto $I_s$, so
we can make choices so that $a_{1,t-1}$, $a_{2,t}$, $a_{1,k}$ and
$a_{k+1,t}\in I_s$, for $2\leq k\leq {t-2}$. This implies that
$\sum\limits_{k=1}^{t-1} {\bar a}_{1,k}a_{k+1,t}$ is a closed
element in the ideal $I_s$ and hence it is exact since $(\bV,d)$
is $s$--formal.
\end{proof}

Other properties of $s$--formal manifolds are given in the
following lemmas.

\begin{lemma}\label{formal}
Let $M$ be a differentiable manifold of dimension $m$. Then $M$
is formal if and only if $M$ is $m$--formal.
\end{lemma}

\begin{proof}
{}From Theorem~\ref{criterio1} and Definition~\ref{primera} it
follows that if $M$ is formal then is $m$--formal because $M$ is
$s$--formal for all $s$.

Conversely, let us suppose that the differentiable manifold $M$ is
$m$--formal and let $(\bigwedge V,d)$ be an $m$--formal minimal
model of $M$. Because $V$ is a graded vector space, we can define
the spaces $N^i$ by $N^i$ = $V^i$ for $i>m$. Denote by $N$ the
graded space $N = \bigoplus\limits_{j>0} N^j$. To prove that $M$
is formal we use Theorem~\ref{criterio1}. It is sufficient to show
that any closed element in the ideal $I(N)$ generated by $N$ in
$\bigwedge V$ is exact. Let $\eta$ be such an element. There are
two possibilities according to $\deg(\eta)\leq m$ or
$\deg(\eta)>m$. If $\deg(\eta)\leq m$, then $\eta$ lies in the
ideal $I_m(N^{\leq m})$, and so $\eta$ is exact because $M$ is
$m$--formal. If $\deg(\eta)>m$, then $\eta$ defines a cohomology
class in the cohomology group $H^{\deg (\eta)}(\bigwedge V)$.
Since that $M$ has dimension $m<\deg (\eta)$ and $(\bigwedge V,d)$
is a model of $M$, the group $H^{\deg(\eta)}(\bigwedge V)$ must be
equal to zero, which implies that the cohomology class $[\eta]$ is
the trivial class and so $\eta$ is exact.
\end{proof}

\begin{lemma}\label{producto}
Let $M_1$ and $M_2$ be differentiable manifolds. For any $s\geq
0$, the product manifold $M=M_1\times M_2$ is $s$--formal if and
only if $M_1$ and $M_2$ are $s$--formal.
\end{lemma}

\begin{proof}
Denote by $(\bigwedge V_i,d_i)$ the minimal model of $M_i$. Then the
minimal model of $M$ is $(\bigwedge V,d)$ with $\bigwedge V=
\bigwedge V_1 \otimes \bigwedge V_2$ and differential $d=d_1 \otimes
1 + 1 \otimes d_2$. Since $M_i$ ($i=1,2$) is $s$--formal,
Lemma~\ref{topologia} implies the existence of a map of differential
algebras
 $$
 \varphi_i\colon \bV_i^{\leq s} \too H^*(M_i)\, ,
 $$
such that the induced map in cohomology $\varphi_i^*$ equals the map
induced by the inclusion $\bV_i^{\leq s} \inc \bV_i$. Consider
$\varphi=\varphi_1\otimes \varphi_2$. Then
 $$
 \varphi \colon \bigwedge (V_1^{\leq s}\oplus V_2^{\leq s}) \too
 H^*(M_1)\otimes
 H^*(M_2)=H^*(M)
 $$
satisfies the conditions of Lemma~\ref{topologia}, hence $M$ is
$s$--formal.

Suppose now that $M=M_1\times M_2$ is $s$--formal. Then by
Lemma~\ref{topologia}, there exists a map of differential algebras
  $$
  \varphi \colon \bigwedge V^{\leq s} \too H^*(M),
  $$
such that $\varphi^*\colon H^*(\bV^{\leq s}) \to H^*(\bV)\cong
H^*(M)$ is the map induced by the inclusion $\bV^{\leq s} \inc \bV$.

Define $\varphi_1$ as the inclusion $\bV_1^{\leq s}\inc \bV^{\leq
s}$ followed by $\varphi$ and by the projection
$H^*(M)=H^*(M_1)\otimes H^*(M_2) \surj H^*(M_1)$. This is a map of
differential algebras and it is easy to see that it satisfies the
conditions of Lemma~\ref{topologia}. In fact, $\varphi_1^*$ equals
the map induced in cohomology by the composition $\bV_1^{\leq s}\inc
\bV^{\leq s} \inc \bV = \bigwedge (V_1\oplus V_2) \surj \bV_1$, and
this map is the inclusion $\bV_1^{\leq s}\inc \bV_1$.
\end{proof}

\section{Formality and $s$--formality}\label{maintheorem}

The purpose of this section is to prove the following theorem.

\begin{theorem}\label{criterio2}
Let $M$ be a connected and orientable compact differentiable
manifold of dimension $2n$, or $(2n-1)$. Then $M$ is formal if and
only if is $(n-1)$--formal.
\end{theorem}

\begin{proof}
In one direction the proof is obvious. So we need to show that the
$(n-1)$--formality of $M$ implies its formality. First suppose
that the theorem holds for any $(n-1)$--formal manifold of
dimension $2n$. Now, if M is a $(n-1)$--formal manifold of $\dim
M= (2n-1)$, the product manifold $M\times S^1$ is $2n$--dimensional
and $(n-1)$--formal according to Lemma~\ref{producto}. Our
assumption implies that $M\times S^1$ is formal. But a product
manifold $M_1\times M_2$ is formal if and only if each one of the
manifolds $M_1$ and $M_2$ is formal. Therefore, $M$ must be
formal, which proves the theorem for odd-dimensional
differentiable manifolds.

To prove the theorem when $\dim M = 2n$ we will show that $M$ is
$(n+r-1)$--formal for any $r\geq 0$ proceeding by induction on
$r$. If $r=0$ then $M$ is $(n-1)$--formal by the hypothesis of the
theorem. Let us suppose that $M$ is $(n+r-1)$--formal and we will
show that $M$ is $(n+r)$--formal for $r\geq 0$.

Let $(\bigwedge V,d)$ be a $(n+r-1)$--formal minimal model of $M$.
By the induction hypothesis, we know that each one of the spaces
$V^{\leq (n+r-1)}$, of generators of degree $\leq (n+r-1)$,
satisfies the conditions of Definition~\ref{primera}.
Since ($\bV,d$) is a minimal differential algebra,
it is possible to order the
generators $\{x_1,x_2,\ldots\}$
of $V^{n+r}$ in such way that
$d\,x_j\in \bigwedge (V^{\leq (n+r-1)}\oplus \la
x_1,\ldots,x_{j-1}\ra)$ for $j\geq 1$. Now, for each generator
$x_i$ of $V^{n+r}$, we define the space $V_i$ by
 $$
 V_i = V^{\leq (n+r-1)}\oplus \la x_1,\ldots,x_i\ra.
 $$
Here we can take $i\geq 0$ and $V_0 = V^{\leq (n+r-1)}$.

We aim to construct a $(n+r)$--formal minimal model of $M$. For
this, for each $x_i\in V^{n+r}$ we shall find $\psi_i\in \bV_{i-1}$
such that
$\hat{x}_i=x_i-\psi_i$ gives a new set of generators, and the space
 $$
 \hV_i = V^{\leq (n+r-1)}\oplus \la \hat{x}_1,\ldots,\hat{x}_i\ra
 $$
satisfies the conditions of Definition~\ref{primera}, i.e.,
$\hV_i$ decomposes as a direct sum $\hV_i = \hC_i\oplus \hN_i$
with $d(\hC_i) = 0$, $d$ injective on $\hN_i$, and that every
closed element in the ideal $I(\hN_i)$, generated by $\hN_i$ in
$\bigwedge \hV_i$, is exact in $\bigwedge V$. Note that $\bigwedge
\hV_i=\bV_i$ for all $i$. If we do this, then $\hV^{n+r}=\la
\hat{x}_1, \hat{x}_2, \ldots \ra$ satisfies the conditions of
Definition~\ref{primera}, and thus $\bV^{\leq (n+r)}=\bigwedge
(V^{\leq (n+r-1)} \oplus \hV^{n+r})$ is $(n+r)$--formal.
We shall proceed by induction on $i$. It is clear for $i=0$. Let us
suppose that it is true for $i-1$, and we shall show it for $i$.

To start with, consider the composition
 \begin{equation} \label{eqn:v1}
 V^{n+r} \stackrel{d}{\too} \bigwedge V^{\leq (n+r)} \too \frac{\bigwedge
 V^{\leq (n+r)}}{d( \bV^{\leq (n+r-1)})}.
 \end{equation}

We reorder the generators of $V^{n+r}$ as follows. Let
$x_1,\ldots, x_p$ be generators of the kernel of~(\ref{eqn:v1}).
Without loss of generality, we may suppose that they are the first
$p$ generators of $V^{n+r}$. Then for $x_i$, $1\leq i \leq p$, we
have that $dx_i$ lies in $d(\bV^{\leq (n+r-1)})$, i.e., there is
some $\psi_i \in \bV^{\leq (n+r-1)}$ with $dx_i=d\psi_i$. Put
$\hat{x}_i=x_i -\psi_i$, so that $d\hat{x}_i=0$.

For $1\leq i \leq p$, define
 $$
 \hC_i= \hC_{i-1}\oplus \la \hat{x}_i\ra=C^{\leq (n+r-1)} \oplus
 \la \hat{x}_1,\ldots, \hat{x}_i\ra
 $$
and
 $$
 \hN_i = \hN_{i-1}=N^{\leq (n+r-1)}.
 $$
Then, we only must show that any closed element in the ideal
$\hN_{i}\cdot ( \bigwedge \hV_i)= \hN_{i-1}\cdot \bigwedge
(\hV_{i-1} \oplus \la \hat{x}_i \ra )$ is exact. Let $\eta\in
\hN_{i}\cdot(\bigwedge \hV_i)$ be a closed element. Thus, $\eta =
\eta_0 + \eta_1 \cdot \hat{x}_i+ \cdots+\eta_k \cdot
{\hat{x}_i}^k$ for some $\eta_j\in \hN_{i-1} \cdot (\bigwedge
\hV_{i-1})$. Moreover, $d\,\eta = 0$ implies that $d\eta_0 +
d\eta_1 \cdot \hat{x}_i+\cdots+d\eta_k \cdot {\hat{x}_i}^k = 0$ in
$\bigwedge \hV_i$. Therefore $d\eta_j = 0$ for $0\leq j\leq k$.
{}From this fact and the induction hypothesis on $\hV_{i-1}$ it
follows that each element $\eta_j$ is exact, and so $\eta$ also is
exact.

\medskip

Now let $i>p$. Then we put $\hC_i=\hC_{i-1}$. We want to see that
there is an element $\psi_i\in \bigwedge\hV_{i-1}$ such that
putting $\hat{x}_i=x_i-\psi_i$ and $\hN_i=\hN_{i-1} \oplus \la
\hat{x}_i\ra$, the decomposition $\hV_i=\hC_i\oplus \hN_i$
satisfies the conditions of Definition \ref{primera}. No matter
$\psi_i$, $d$ is injective in $\hN_i$. This follows from the fact
that~(\ref{eqn:v1}) is injective in $\la x_{p+1},\cdots, x_i\ra$
and that $\hat{x}_j=x_j-\psi_j$ with $\psi_j\in \bigwedge
\hV_{j-1}$, for $p+1\leq j\leq i$.

For the time being, write $N_i=\hN_{i-1}\oplus \la x_i\ra$
and 
let $\eta\in N_i \cdot \bigwedge (\hV_{i-1}\oplus \la x_i \ra)=
(\hN_{i-1} \oplus \la x_i \ra) \cdot {\bigwedge (\hV_{i-1} \oplus
\la x_i\ra)}$ be a closed element. Then $\eta = \eta_0 + \eta_1
\cdot x_i+\cdots+\eta_k \cdot {x_i}^k$. We distinguish three
cases:

\begin{enumerate}
\item $k=0$. Now $\eta = \eta_0$ with
$\eta_0\in \hN_{i-1} \cdot (\bigwedge \hV_{i-1})$. By the
induction hypothesis on $i$ we know that $\eta$ is exact.

\item $k\geq 2$. Note that in this case the degree of $x_i$ must
be even. Because $\deg\,x_i\geq n$ we have that $\deg {\eta}\geq
{2n}$. If either $k>2$ or $k=2$ and $\deg\,x_i>n$, it happens that
$\deg {\eta}>2n$. Then $\eta$ must be exact because $H^{>2n}(M)=0$
and $(\bV, d)$ is a minimal model of $M$. The only remaining
possibility is that $k=2$ and $\deg\,x_i = n$. In this case $\deg
{\eta}=2n$ and $\eta$ has an expression of the form $\eta = \eta_0
+ \eta_1 \cdot x_i+\lambda\, {x_i}^2$ where $\lambda$ is a
non-zero real number. Thus $0=d\eta =(d\eta_0 +\eta_1 \cdot dx_i)+
(d\eta_1 +2\lambda\, dx_i) x_i$ and hence $\eta_1 + 2\lambda\,
x_i$ is closed.

Now $\eta_1\in (\bigwedge \hV_{i-1})^n$, so it must be of the form
$\eta_1=a+b$ with $a\in \bV^{\leq (n-1)}$, $b\in \la {x}_1,\ldots,
{x}_{i-1}\ra$. So $d(2\lambda\, x_i- b)=-da$. This means that
$2\lambda\, x_i- b$ is in the kernel of the map~(\ref{eqn:v1}),
which is a contradiction.

\item $k=1$. Thus $\eta = \eta_0+ \eta_1 \cdot x_i$,
with $\eta_0\in \hN_{i-1} \cdot (\bigwedge \hV_{i-1})$ and $\eta_1
\in \bigwedge \hV_{i-1}$. In this case $d\eta=0$ implies that
$d\eta_1=0$. We shall see that we can change $x_i$ to some
$\hat{x}_i=x_i-\psi_i$ with an element $\psi_i\in (\bigwedge
\hV_{i-1})^{n+r}$ so that any closed element of the form
$\eta_0+\eta_1 \cdot \hat{x}_i$ must be exact in $\bigwedge V$.
Note that substituting $x_i$ by $\hat{x}_i$ does not spoil the
argument in the previous two cases.

If $\deg {\eta}>2n$ and $\eta$ is closed, one has that $\eta$ is
exact by the same argument as the case $\deg {\eta}>2n$ of (ii).

Now, we deal with the case that $\eta$ has degree $2n$. This
implies that $\eta_1$ is closed of degree $n-r$. 
In order to show the exactness of $\eta$ we proceed as follows. In
general,  consider the collection of those closed $z_j \in
(\bigwedge \hV_{i-1})^{n-r}$ such that there exists $\kappa_j \in
\left( \hN_{i-1}\cdot(\bigwedge \hV_{i-1}) \right)^{2n}$ in such
way that the element $z_j \cdot x_i +\kappa_j$ is closed. Hence
there is $\xi_j \in \bV$ satisfying
  \begin{equation} \label{eqn:v2}
  z_j \cdot x_i +\kappa_j = \lambda_j \omega + d\xi_j
  \end{equation}
where $\lambda_j$ are real numbers and $\omega$ is a (fixed)
closed element of degree $2n$ generating $H^{2n}(\bigwedge V)\cong
\RR$. We want to achieve that $\lambda_j=0$ for all $j$. First,
for a given $z_j$, suppose that we have two different expressions
$z_j \cdot x_i +\kappa_j = \lambda_j \omega + d\xi_j$ and $z_j
\cdot x_i +\kappa_j' = \lambda_j' \omega + d\xi_j'$. Then the
difference $\kappa_j-\kappa_j'=(\lambda_j -\lambda_j')\omega +
d\xi_j-d\xi_j'$ is closed and lives in $\hN_{i-1} \cdot (\bigwedge
\hV_{i-1})$. By induction hypothesis, it is exact and hence
$\lambda_j=\lambda_j'$. So if we manage to make $\lambda_j=0$ we
have dealt with any possible expression (\ref{eqn:v2}) for $z_j$.

So we may restrict to a basis of those $z_j$ satisfying
(\ref{eqn:v2}). If $[z_j]=0$ (for example, when $H^{n-r}(\bV)=0$)
then $z_j=d \phi$, with $\phi\in (\bigwedge V)^{n-r-1}$. Clearly
one can take $\phi \in N^{\leq(n-r-1)}\cdot \bigwedge V^{\leq
(n-r-1)}$. Now
 $$
 z_j \cdot  x_i +\kappa_j= d\phi \cdot  x_i +\kappa_j=
 d(\phi \cdot  x_i) - (-1)^{n-r-1} \phi \cdot  dx_i +\kappa_j,
 $$
which implies that $\phi \cdot  dx_i + (-1)^{n-r}\kappa_j \in
\hN_{i-1}\cdot (\bigwedge \hV_{i-1})$ is closed and hence exact
taking into account (\ref{eqn:v2}). So $z_j \cdot x_i+\kappa_j$ is
exact and $\lambda_j=0$. This means that if $\eta = \eta_0+ \eta_1
\cdot x_i$ is closed and $[\eta_1]=0$ then $\eta$ is exact.

Therefore we may restrict ourselves to a collection of $z_j$ such
that $[z_j]$ are a basis of the possible $z_j$'s
satisfying~(\ref{eqn:v2}).

Let $[z_1],\ldots, [z_k] \in H^{n-r}(\bigwedge V)$ be such a
basis. By Poincar\'{e} duality there is some element $[\psi_i] \in
H^{n+r}(\bigwedge V)$ with $[z_j] \cdot [\psi_i]=\lambda_j
[\omega]$ for all $j$. Such a closed element $\psi_i\in (\bigwedge
V)^{n+r}$ must lie in $\bigwedge ( V^{<(n+r)}\oplus \la
\hat{x}_{1},\ldots,\hat{x}_{p}\ra)$
since $(\bigwedge V)^{n+r}=(\bigwedge V^{<(n+r)})^{n+r}\oplus \la
\hat{x}_{1},\ldots,\hat{x}_{p}\ra \oplus \la
x_{p+1},x_{p+2},\ldots\ra$ and, by~(\ref{eqn:v1}), $\psi_i$ cannot
have a component in $\la x_{p+1},x_{p+2},\ldots\ra$. 
Now define $\hat{x}_i=x_i-\psi_i$. Then, from (\ref{eqn:v2}), it
is easy to check that $z_j \cdot \hat{x}_i+\kappa_j$ is exact,
i.e., whenever
 $$
 z \cdot  \hat{x}_i+\kappa, \quad z \hbox{ closed in }
 (\bigwedge \hV_{i-1})^{n-r},
 \quad  \kappa \in \left( \hN_{i-1}\cdot (\bigwedge \hV_{i-1})
 \right)^{2n},
 $$
is closed, it is exact.

Put $\hN_i=\hN_{i-1}\oplus \la \hat{x}_i\ra$. Now we can check
that the conditions of Definition \ref{primera} hold. It only
remains to show that if $\eta = \eta_0+ \eta_1 \cdot \hat{x}_i$ is
closed of degree $<2n$, with $\eta_0\in \hN_{i-1}\cdot (\bigwedge
\hV_{i-1})$ and $\eta_1 \in \bigwedge \hV_{i-1}$, then $\eta$ is
exact. But $[\eta_1]\in H^{p}(\bigwedge V)$ with $p<n-r$. If
$H^{n-r-p}(\bigwedge V)=0$, then $H^{n+r+p}(\bigwedge V)=0$ and
thus $[\eta_0+\eta_1 \cdot  \hat{x}_i]=0$. If $H^{n-r-p}(\bigwedge
V) \not= 0$, take an arbitrary $[w]\in H^{n-r-p}(\bigwedge V)$.
Now we have
 $$
  [w]\cdot [\eta_0+\eta_1 \cdot  \hat{x}_i] = [w \cdot  \eta_0 + w \cdot
  \eta_1 \cdot  \hat{x}_i],
 $$
with $w \cdot  \eta_0\in \hN_{i-1}\cdot (\bigwedge \hV_{i-1})$ and
$w \cdot  \eta_1$ closed of degree $n-r$. Then by the construction
above, $[w]\cdot [\eta_0+\eta_1 \cdot \hat{x}_i] =0$. Since this
holds for arbitrary $[w]$, Poincar\'{e} duality implies that
$[\eta_0+\eta_1 \cdot \hat{x}_i] =0$.
\end{enumerate}

\end{proof}

Miller's theorem~\cite{Mi} for the formality of a
$(k-1)$--connected compact manifold of dimension less than or
equal to $(4k-2)$ follows easily from our Theorem~\ref{criterio2}.

\begin{theorem}~\cite{Mi} \label{Miller}
Let $M$ be a $(k-1)$--connected compact manifold of dimension less
than or equal to $(4k-2)$, $k>1$. Then $M$ is formal.
\end{theorem}

\begin{proof}
Since $M$ is $(k-1)$--connected, a minimal model ($\bigwedge V,d)$
of $M$ must satisfy $V^i = 0$ for $i\leq {k-1}$ and $V^k = C^k$
(i.e.,\ $N^k=0$). Therefore the first non-zero differential, being
decomposable, must be $d: V^{2k-1} \to V^k \cdot V^k$. This
implies that $V^j = C^j$ (i.e.,\ $N^j=0$) for $k\leq j\leq
{(2k-2)}$. Hence $M$ is $2(k-1)$--formal. Now, using
Theorem~\ref{criterio2} we have that $M$ is formal. Note that $M$
is orientable since it is simply connected.
\end{proof}

Note that this implies in particular that any simply connected
compact manifold of dimension less than or equal to 6 is formal,
which is a result of~\cite{NM} previous to Miller's theorem.

Also, as a consequence of Theorem~\ref{criterio2} we have the corollary
following.

\begin{corollary}\label{8-manifolds}
Any simply connected compact manifold, of arbitrary dimension, is
$2$--formal. Moreover, a simply connected compact manifold $M$ of
dimension $7$ or $8$ is formal if and only if it is $3$--formal.
\end{corollary}


\begin{remark}
Theorem \ref{criterio2} continues to hold for rational Poincar\'{e}
duality spaces (see Remark~\ref{rem:cw}).
\end{remark}

A symplectic manifold $(M,\omega)$ is said to be {\it
symplectically aspherical\/} if $\omega|_{\pi_2(M)}=0$, that is,
 $$
 \int_{S^2} f^*\omega = 0
 $$
for every map $f\colon S^2 \to M$. Hurewicz's theorem implies that
a compact symplectically aspherical manifold always has a
non-trivial fundamental group.

\begin{remark}
We note that, from Theorem~\ref{criterio2}, if $M$ is a
$1$--formal manifold of dimension $4$, then $M$ is formal, so all
Massey products are trivial. However, the converse is not true.
Actually, Amor\'os and Kotschick~\cite{AK} claim to have examples of
non-formal manifolds of dimension $4$ with all Massey products of length
$t\leq K$ vanishing, for any arbitrary large number $K$, as well as to
have examples of non-formal symplectic $4$--manifolds which are
hard Lefschetz. Their examples are symplectically aspherical (and
hence non-simply connected).
\end{remark}

\section{Lefschetz property} \label{sectionlefschetz}

In this section we introduce the $s$--Lefschetz property for any
compact symplectic manifold, generalizing the hard Lefschetz
property. We will study this property for Donaldson submanifolds
of  symplectic manifolds in the next section.

\begin{definition} \label{s-Lefschetz}
  Let $(M,\omega)$ be a compact symplectic manifold of dimension $2n$. We
  say that $M$ is $s$--Lefschetz \textup{}with $s\leq (n-1)$\textup{} if
  $$
   [\omega]^{n-i}: H^i(M) \too H^{2n-i}(M)
  $$
  is an isomorphism for all $i\leq s$. By extension, if we say that
  $M$ is $s$--Lefschetz with $s\geq n$ then we just mean that $M$
  is hard Lefschetz.
\end{definition}

Note that $M$ is $(n-1)$--Lefschetz if $M$ satisfies the hard
Lefschetz theorem. Also it is said in~\cite{McD2} that $M$ is a
Lefschetz manifold meaning that $M$ is $1$--Lefschetz. $M$ is
$0$--Lefschetz if it is cohomologically symplectic.

In Section~\ref{examples} we present examples of compact
symplectic manifolds which are $s$--Lefschetz but not
$(s+1)$--Lefschetz for $s=0,1,2$ (see Examples 1, 3 and 4).
However, we do not know examples of symplectic manifolds being
$s$--Lefschetz but not $(s+1)$--Lefschetz for $s\geq 3$.

\begin{proposition} \label{prod-Lefschetz}
 Let $(M_1, \omega_1)$ and $(M_2,\omega_2)$ be symplectic
 manifolds and let
 $M=M_1\times M_2$ with the symplectic form $\omega=\lambda_1 \omega_1
 + \lambda_2 \omega_2$, $\lambda_1,\lambda_2$ non-zero real
 numbers. Then for any $s$, $M$ is $s$--Lefschetz if and
 only if $M_1$ and $M_2$ are $s$--Lefschetz.
\end{proposition}

\begin{proof}
 First note that we may rescale the symplectic forms $\omega_1$
 and $\omega_2$ as $\lambda_1\omega_1$ and $\lambda_2\omega_2$, so
 the coefficients can be supposed equal to one. Also let $2n_1$
 and $2n_2$ be the dimensions of $M_1$ and $M_2$ respectively. Let
 $n=n_1+n_2$.

 Suppose first that $H^*(M)$ is $s$--Lefschetz. Let us see that $M_1$ is
 also $s$--Lefschetz. Take $i\leq s$ with $i\leq n_1-1$. Let
 $a_i\in H^i(M_1)$ and
 suppose that $[\omega_1]^{n_1-i}a_i=0$. It is enough to see that
 $a_i=0$ because then $[\omega_1]^{n_1-i}\colon H^i(M_1)\to
 H^{n_1-i}(M_1)$ is injective and hence an isomorphism. But
 $$
 [\omega]^{n-i} (a_i\otimes 1) = \sum_{k\geq 0} {n-i \choose n_2-k}
 [\omega_1]^{n_1-i+k} a_i \otimes [\omega_2]^{n_2-k} =0
 $$
 and the $s$--Lefschetz property for $M$ implies that $a_1\otimes
 1=0$ and hence $a_i=0$.

 For the converse, the $s$--Lefschetz property for $M_1$ implies that we
 may decompose $H^*(M_1)=(\oplus P_i) \oplus R_1$ in vector subspaces,
 so that
 $$
 P_i=\la e_{i}, [\omega_1] e_{i}, \ldots, [\omega_1]^{n_1-d_i}
 e_{i}\ra,
 $$
 where $e_i\in H^{d_i}(M_1)$, $d_i\leq s$, and $[\omega_1]^{n_1-d_i+1} e_{i}=0$.
 This is possible thanks to
 the $s$--Lefschetz property. The elements $e_{i}$ are called
 primitive elements. The subspace $R_1$ is
 concentrated in degrees going from $s+1$ up to $2n_1-s-1$.
 Similarly $H^*(M_2)=(\oplus Q_j) \oplus R_2$, where
 $Q_j=\la f_j, [\omega_2] f_j, \ldots, [\omega_2]^{n_2-d_j}
 f_j\ra$ and $R_2$ is
 concentrated in degrees going from $s+1$ up to $2n_2-s-1$.
 Therefore
 $$
 H^*(M_1\times M_2)= \left(\bigoplus_{i,j} P_i\otimes Q_j \right) \oplus R,
 $$
 where $R=R_1\otimes H^*(M_2) + H^*(M_1) \otimes R_2$ (not a
 direct sum). Then $R$ is concentrated in degrees going from $s+1$
 up to $2n_1+2n_2 -s-1$. This means that $R$ is irrelevant for
 checking the $s$--Lefschetz condition for $M$.

 On the other hand,
 $$
 P_i \otimes Q_j=\la e_i \otimes f_j,\ldots,  [\omega_1]^a e_i \otimes
 [\omega_2]^b f_j, \ldots, [\omega_1]^{n_1-d_i}e_i\otimes
 [\omega_2]^{n_2-d_j} f_j \ra
 $$
 satisfies the hard Lefschetz condition with respect to
 $\omega=\omega_1+\omega_2$. Therefore $M$ is $s$--Lefschetz.
\end{proof}

\section{Donaldson submanifolds of  symplectic manifolds}
\label{sectionsubvarities}
In this section we study the conditions under which Donaldson
symplectic submanifolds are formal and/or satisfy the hard
Lefschetz theorem.

Let $(M,\omega)$ be a compact symplectic manifold of dimension
$2n$ with $[\omega]\in H^2(M)$ admitting a lift to an integral
cohomology class. In \cite{Do} Donaldson constructs symplectic
submanifolds $Z\inc M$ of dimension $2n-2$ whose Poincar\'e dual
$\PD[Z]=k[\omega]$ for any large multiple of $[\omega]$. Moreover,
these submanifolds satisfy a {\em Lefschetz theorem in hyperplane
sections\/}, meaning that the inclusion $j\colon Z\inc M$ is
$(n-1)$--connected. In particular, the map there $j^*\colon
H^i(M)\to H^i(Z)$ is an isomorphism for $i<n-1$ and a monomorphism
for $i=n-1$.

More in general, let $X$ and $Y$ be compact manifolds. We say that
a differentiable map $f\colon X \to Y$ is a {\it homology
$s$--equivalence} ($s\geq 0$) if it induces isomorphisms
$f^*\colon H^i(Y) \isom H^i(X)$ on cohomology for $i<s$, and a
monomorphism $f^*\colon H^s(Y)\inc H^s(X)$ for $i=s$. Therefore
$Z\inc M$ is a homology $(n-1)$--equivalence.

In~\cite{DGMS} it is proved that if $F\colon B_1 \to B_2$ is a
morphism of differential algebras inducing an isomorphism on
cohomology, and $\rho_i \colon A_i \to B_i$ is a minimal model for
$B_i$ ($i=1,2)$, then $F$ induces $\hat{F} \colon A_1
\longrightarrow A_2$, unique up to homotopy, subject to the
condition $F\circ \rho_1 = \rho_2 \circ \hat{F}$. For a homology
$s$--equivalence we have:

\begin{proposition}\label{equivalence1}
Let $X$ and $Y$ be compact manifolds and let $f\colon X \to Y$ be
a homology  $s$--equivalence. Then there exist minimal models
$(\bV_X,d)$ and $(\bV_Y,d)$ of $X$ and $Y$, respectively, such
that $f$ induces a morphism of differential algebras $F\colon
(\bigwedge V_Y^{\leq s},d) \to (\bigwedge V_X^{\leq s},d)$ where
$F: V_{Y}^{<s} \isom V_{X}^{<s}$ is an isomorphism and $F: V_{Y}^s
\inc V_{X}^s$ is a monomorphism.
\end{proposition}

\begin{proof}
We can do this by induction on $s$, being evident for $s=0$. So we
can suppose that if $f$ is a homology $(s-1)$--equivalence, there
exist minimal models $(\bV_{X},d)$ and $(\bV_{Y},d)$ for $X$ and
$Y$, respectively, and a morphism $F \colon (\bigwedge V_Y^{<s},d)
\longrightarrow (\bigwedge V_X^{<s},d)$ such that $F\colon
V_{Y}^{<(s-1)} \isom V_{X}^{<(s-1)}$ is an isomorphism and
$F\colon V_{Y}^{s-1} \inc V_{X}^{s-1}$ is a monomorphism. Then, we
shall prove the Proposition for a homology $s$--equivalence $f$.
So, $f^*$ induces: $H^{s-1}(Y) \iso H^{s-1}(X)$ and $H^s(Y) \inc
H^s(X)$. We have $H^{j}(Y) \iso H^{j}(\bigwedge V_{Y})$ and
$H^{j}(X) \iso H^{j}(\bigwedge V_{X})$ for any $j$. Hence $f^*$
induces: $H^{s-1}(\bigwedge V_{Y}) \iso H^{s-1}(\bigwedge V_{X})$
and $H^{s}(\bigwedge V_{Y}) \inc H^s(\bigwedge V_{X})$.

For convenience, we shall denote by $\hat{\mu}$ the element
$\hat{\mu} = F (\mu)$, for $\mu \in \bV_{Y}^{\leq (s-1)}$. First,
we prove that $F\colon V_{Y}^{s-1} \isom V_{X}^{s-1}$ is an
isomorphism. Let us order the generators of $V_{X}^{s-1}$ so that
  $$
  V_X^{s-1}=F(V_Y^{s-1}) \oplus \la x_1, x_2,\ldots \ra,
  $$
where $dx_i \in \bigwedge (F(V_Y^{\leq (s-1)}) \oplus \la
x_1,\ldots, x_{i-1}\ra)$. We shall show that there is no second
summand above by showing that $x_1 \in F(V_Y^{\leq (s-1)})$.

We have $dx_1 \in \bigwedge (F(V_Y^{\leq (s-1)}))$. Set
$\hat{\eta} =dx_1=F(\eta)$ for some $\eta \in (\bV_Y^{\leq
(s-1)})^s$. So $d\eta=0$ and hence $[\eta]\in H^{s}(\bV_{Y})$.
Under the monomorphism $f^*\colon H^{s}(\bV_{Y})\inc
H^{s}(\bV_{X})$, we get that $f^*[\eta]= [\hat{\eta}]$. Since
$[\hat{\eta}]=0$ in $H^s(\bV_{X})$, we have $[\eta]=0$ in
$H^s(\bV_{Y})$. This guarantees the existence of some $\alpha \in
(\bigwedge V_Y)^{s-1}\inc (\bigwedge V_X)^{s-1}$ such that $\eta
=d\alpha$. Therefore $d(x_1-\hat{\alpha})=0$ and so there is a
well-defined cohomology class $[x_1-\hat{\alpha}]$ living in
$H^{s-1}(\bV_{X})$. But this space is isomorphic by $f^*$ to
$H^{s-1}(\bV_{Y})$. Hence there exist a cohomology class $[\mu]
\in H^{s-1}(\bV_{Y})$ such that $[\hat{\mu}]=[x_1-\hat{\alpha}]$.
This implies that there are a closed element $\mu \in (\bigwedge
V_Y)^{s-1}$ and $\xi \in (\bigwedge V_Y)^{s-2} \iso (\bigwedge
V_X)^{s-2}$ such that
  $$
  x_1-\hat{\alpha}= \hat{\mu} + d\hat{\xi},
  $$
which means that $x_1 \in F(\bigwedge V_{Y}^{\leq (s-1)})$. So it
must be $V_X^{s-1}\iso V_Y^{s-1}$.

Now let us see that $f$ induces a map $F:V_Y^s \inc V_X^s$. Write
$V_Y^s=\la y_1, y_2,\ldots\ra$ with $dy_i \in \bigwedge (V_Y^{\leq
(s-1)} \oplus \la y_1,\ldots, y_{i-1}\ra)$. Now $F(y_i)\in
(\bV_X)^s= (\bV_X^{<s})^s \oplus V_X^s$. Since we already have
that $\bV_X^{<s} \iso \bV_Y^{<s}$, we may modify $y_i$ by adding a
suitable element in $\bV_Y^{<s}$ (and keep on denoting it by
$y_i$) so that $F(y_i)\in V_X^s$. Now we can assume that $\la
y_1,\ldots, y_{r-1}\ra$ injects into $V_X^s$ but
$\hat{y}_r=F(y_r)=0$. Then, on the one hand, we have
$d\hat{y}_r=0$ , thus $[\hat{y}_r]$ is the zero class in
$H^s(\bV_{X})$. On the other hand,
  \begin{equation} \label{eqn:v3}
  dy_r=P(y_1, \ldots, y_{r-1})\in \bigwedge (V_Y^{\leq (s-1)} \oplus
  \la y_1,\ldots, y_{r-1}\ra),
  \end{equation}
that is, $dy_r$ is a polynomial in previous generators.

Applying $F$ to~(\ref{eqn:v3}), we have $0=P(\hat{y}_1,\ldots,
\hat{y}_{r-1})$ in the free algebra $\bigwedge (V_X^{\leq (s-1)}
\oplus \la \hat{y}_1,\ldots,\hat{y}_{r-1}\ra)$. Therefore
$dy_r=0$, so $[y_r]$ is a cohomology class in $H^s(\bV_{Y})$. As
$H^s(\bV_{Y})\inc H^s(\bV_{X})$ and $[\hat{y}_r]=0$ in
$H^s(\bV_{X})$, we get that $[y_r]=0$, i.e., $y_r=d\eta$ for
$\eta\in \bigwedge V_Y^{\leq (s-1)} \iso \bigwedge V_X^{\leq
(s-1)}$. Applying $F$ to $y_r=d\eta$, we have that $0=d\hat{\eta}$
in $\bigwedge V_X$, so $[\hat{\eta}]\in H^{s-1}(\bV_{X})\iso
H^{s-1}(\bV_{Y})$. Therefore there are a closed element $g \in
(\bigwedge V_Y)^{s-1}$ and $\xi \in (\bigwedge V_Y)^{s-2} \iso
(\bigwedge V_X)^{s-2}$ such that
  $$
  \eta=g +d\xi.
  $$
Hence $d\eta=0$ in $\bigwedge V_Y$ and $y_r=0$. This is a
contradiction since $y_r$ is a generator. So it must be that $\la
y_1,\ldots, y_{r}\ra \inc V_X^s$. This implies that $V_Y^s\inc
V_X^s$.
\end{proof}

\begin{theorem}\label{equivalence2}
\begin{enumerate}
\item Let $X$ and $Y$ be compact manifolds, and let $f \colon X
\to Y$ be a homology $s$--equivalence. If $Y$ is $(s-1)$--formal
then $X$ is $(s-1)$--formal. \item Let $M$ be a compact symplectic
manifold of dimension $2n$ and let $Z\inc M$ be a Donaldson
submanifold. For each $s\leq n-2$, if $M$ is $s$--formal then $Z$
is $s$--formal. In particular, $Z$ is formal if $M$ is
$(n-2)$--formal.
\end{enumerate}
\end{theorem}

\begin{proof}
Let $(\bigwedge V_X,d)$ and $(\bigwedge V_Y,d)$ be the minimal
models of $X$ and $Y$, respectively, constructed in
Proposition~\ref{equivalence1}. For $i<s$, decompose $V_Y^i=C_Y^i
\oplus N_Y^i$ satisfying the conditions of Definition
\ref{primera}. Then, taking into account
Proposition~\ref{equivalence1}, we set $V_X^i=C_X^i \oplus N_X^i$
under the natural isomorphism  $F\colon V_Y^i\iso V_X^i$, $i<s$.
Consider a closed element $F(\eta)=\hat{\eta} \in N_X^{<s}\cdot
(\bigwedge V_X^{<s})$. Hence $\eta$ is a closed element in
$N_Y^{<s}\cdot (\bigwedge V_Y^{<s})$ and, by the
$(s-1)$--formality of $Y$, it is exact, i.e., $\eta=d\xi$, for
$\xi\in \bigwedge V_Y$. Take the image $\hat{\eta}=d(F(\xi))$ in
$\bigwedge V_X$. This proves (i). Now (ii) follows from (i) and
using that the inclusion $j\colon Z\inc M$ is a homology
$(n-1)$--equivalence.
\end{proof}

\begin{theorem} \label{Lefschetz-Donaldson}
Let $M$ be a compact symplectic manifold of dimension $2n$, and
let $Z\inc M$ be a Donaldson submanifold. Then, for each $s\leq
n-2$, $M$ is $s$--Lefschetz if and only if $Z$ is $s$--Lefschetz.
\end{theorem}

\begin{proof}
For any differential form $x$ on $M$, we shall denote by $\hat{x}$
the differential form on $Z$ given by $\hat{x}=j^*(x)$. Let
$p=2(n-1)-i$, where $i \leq {(n-2)}$, and consider the restriction map
$j^*\colon H^p(M)\to H^p(Z)$. Then, for $[z]\in H^p(M)$, we claim that
 \begin{equation} \label{eqn:v4}
  j^*[z]=0  \iff  [z]\cup [\omega]=0.
 \end{equation}
This can be seen via Poincar\'{e} duality. Clearly $j^*[z]=0$ if and
only if for any $a\in H^i(Z)$ we have $j^*[z]\cdot a=0$. We know that
there is an isomorphism $H^i(Z)\iso H^i(M)$ ($i\leq n-2$), thus we can
assume that there is a closed $i$--form $x$ on $M$ with
$[x|_Z]=[\hat{x}]=a$. So
 $$
 j^*[z]\cdot [\hat{x}]=\int_Z \hat{z}\wedge \hat{x} =\int_M z\wedge  x \wedge
 k\omega,
 $$
 since $[Z]=k \PD [\omega]$. Hence $j^*[z]=0$ if and only if
 $[z\wedge \omega] \cdot [x]=0$ for all $[x]\in
 H^i(M)$, from where the claim follows.

Now suppose that $M$ is $s$--Lefschetz, so
$[\omega]^{n-i}:H^i(M)\to H^{2n-i}(M)$ is an isomorphism for
$i\leq s$. We want to check that the map
$[\omega_Z]^{n-1-i}:H^i(Z) \to
 H^{2n-2-i}(Z)$ is injective. Let $[\hat{z}]\in H^i(Z) \iso H^i(M)$ and
 extend it to $[z] \in H^i(M)$. Then, $[\omega_Z]^{n-1-i} [\hat{z}]=0$
 implies that $j^*[\omega^{n-1-i} \wedge z]=0$, which by~(\ref{eqn:v4})
 is equivalent to
 $[\omega^{n-1-i} \wedge z\wedge \omega]=0$. Using the
 $s$--Lefschetz property of $M$, we get $[z]=0$ and thus
 $[\hat{z}]=0$.

 The converse is easy. If $Z$ is $s$--Lefschetz and we take $[z]\in
 H^i(M)$ such that $[\omega^{n-i}\wedge z]=0$, from~(\ref{eqn:v4})
 it follows that
 $j^*[\omega^{n-1-i} \wedge z]=0$, i.e., $[\omega_Z^{n-1-i} \wedge
 z|_Z]=0$.  Hence $[\hat{z}]=0$ in $H^i(Z)$ and so $[z]=0$ since
 $i\leq n-2$.
\end{proof}


\begin{corollary} \label{hardLefschetzDonaldson}
Let $M$ be a compact symplectic manifold of dimension $2n$, and
let $Z\inc M$ be a Donaldson submanifold. We have that if $M$ is
hard Lefschetz, $Z$ is also hard Lefschetz. Moreover, $M$ is
$(n-2)$--Lefschetz (but not necessarily hard Lefschetz) if and
only if $Z$ is hard Lefschetz.
\end{corollary}

In the Example 3 of Section~\ref{examples} we exhibit examples of
$6$--dimensional compact symplectic manifolds which are
$1$--Lefschetz but not $2$--Lefschetz (i.e.,\ not hard Lefschetz),
so its Donaldson symplectic submanifolds are $1$--Lefschetz and
thus hard Lefschetz.

\begin{corollary} \label{cohomologyDonaldson}
 Under the conditions of Theorem~\ref{Lefschetz-Donaldson},
for each $p=2(n-1)-i$ with $i \leq s$,
there is an isomorphism
 $$
 H^p(Z) \iso \frac{H^p(M)}{\ker([\omega]:H^p(M)\surj
 H^{p+2}(M))}.
 $$
\end{corollary}

\begin{proof}
{}From~(\ref{eqn:v4}), we know that there is an inclusion
  $$
  \frac{H^p(M)}{\ker([\omega]:H^p(M)\to H^{p+2}(M))}
  \hookrightarrow H^p(Z).
  $$
 Furthermore, the map $H^p(M) \to H^{p+2}(M)$ is surjective since the
 $s$--Lefschetz property guarantees an isomorphism $H^{i}(M) \to H^{p+2}(M)$.
 Computing dimensions, we have $b^p(M)-(b^p(M)-b^{p+2}(M))=
 b^{p+2}(M)=b^i(M)=b^i(Z) =b^p(Z)$, which completes the proof.
\end{proof}

\begin{remark}
We must note that under the conditions of
Theorem~\ref{Lefschetz-Donaldson}, if $M$ has a non-zero Massey
product $\la [\alpha_1],[\alpha_2],\dots,[\alpha_t]\ra \subset
H^{r}(M)$, with $r \leq {(n-1)}$, then it defines a non-zero
Massey product $\la
[\hat{\alpha}_1],[\hat{\alpha}_2],\dots,[\hat{\alpha}_t]\ra$ of
$Z$, where $\hat{\alpha}_i=j^*(\alpha_i)$, $1 \leq i \leq t$. If
$r \geq n$ and the cohomology classes $[\hat{\alpha}_i]$, $1 \leq
i \leq t$, are non-trivial in $Z$ then, from
Corollary~\ref{cohomologyDonaldson}, it follows that $\la
[\hat{\alpha}_1],[\hat{\alpha}_2],\dots,[\hat{\alpha}_t]\ra$ is
non-zero if and only if the cup product of $[\omega]$ by any
representative of  $\la [\alpha_1],[\alpha_2],\dots,[\alpha_t]\ra$
is a non-trivial cohomology class of $M$. However, as we show in
Section~\ref{examples}, it can happen that a Massey product $\la
[\hat{\alpha}_1],[\hat{\alpha}_2],\dots,[\hat{\alpha}_t]\ra$ is
defined on $Z$ but $\la [\alpha_1],[\alpha_2],\dots,[\alpha_t]\ra$
is not defined on $M$.
\end{remark}


\section{Examples} \label{examples}

We shall apply the previous results to study the $s$--formality
and the $s$--Lefschetz property of some compact symplectic
manifolds and their Donaldson submanifolds. Five examples will be
developed.

The first one is the well known Kodaira--Thurston manifold $KT$;
it is the simplest nontrivial example of a compact symplectic
manifold with no K\"ahler metric. Example $2$ is the Iwasawa
manifold $I_3$, any Donaldson symplectic submanifold $Z$ of $I_3$
is neither formal nor hard Lefschetz; moreover $Z \inc I_3$ has no
complex structures. Example $4$ allows us to show a Donaldson
symplectic submanifold of dimension eight which is formal simply
connected but not hard Lefschetz. Example $3$ is a
$6$--dimensional compact symplectically aspherical manifold $M$
which is $1$--formal but not $2$--formal, it has the
$1$--Lefschetz property but not the $2$--Lefschetz property; any
Donaldson symplectic submanifold of $M$ is formal and hard
Lefschetz but it does not carry K\"ahler metrics. In Example 5 we
show that there are compact oriented smooth manifolds which are
$s$--formal but not $(s+1)$--formal, for any $s\geq 2$.

{\bf Example 1} {\it The Kodaira--Thurston manifold\/}. Let $H$ be
the Heisenberg group, that is, the connected nilpotent Lie group
of dimension $3$ consisting of matrices of the form
 $$
 a=\pmatrix{1&x&z\cr 0&1&y\cr 0&0&1 \cr},
 $$
where $x,y,z \in {\RR}$. Then a global system of coordinates
${x,y,z}$ for $H$ is given by $x(a)=x$, $y(a)=y$, $z(a)=z$, and a
standard calculation shows that a basis for the left invariant
$1$--forms on $H$ consists of
    $$
    \{dx, dy, dz-xdy\}.
    $$
Let $\Gamma$ be the discrete subgroup of $H$ consisting of
matrices whose entries are integer numbers. So the quotient space
$M^3 =\Gamma{\backslash} H$ is compact. Hence the forms $dx$,
$dy$, $dz-xdy$ descend to $1$--forms $\alpha$, $\beta$, $\gamma$
on $M^3$.

The Kodaira--Thurston manifold $KT$ is the product $KT = M^3
\times S^1$. Then, there are $1$--forms ${\alpha, \beta, \gamma,
\eta}$ on $KT$ such that
 $$
 d\alpha=d\beta=d\eta=0, \quad  d\gamma=-\alpha \wedge \beta,
 $$
and such that at each point of $KT$, $\{\alpha, \beta, \gamma,
\eta\}$ is a basis for the $1$--forms on $KT$. Moreover, it is
easy to use Nomizu's theorem~\cite{No} to compute the real
cohomology of $KT$
\begin{eqnarray*}
 H^0(KT) &=& \la 1\ra, \\
 H^1(KT) &=& \la [\alpha], [\beta], [\eta]\ra,\\
 H^2(KT) &=& \la [\alpha \wedge \gamma], [\beta\wedge \gamma],
   [\alpha \wedge \eta], [\beta \wedge \eta]\ra,\\
 H^3(KT) &=& \la [\alpha \wedge \gamma \wedge \eta],
   [\beta\wedge \gamma \wedge \eta], [\alpha \wedge \beta \wedge \gamma]\ra,\\
 H^4(KT) &=& \la [\alpha \wedge \beta \wedge \gamma \wedge \eta]\ra.
\end{eqnarray*}

Using again Nomizu's theorem, the minimal model of $KT$ is the
differential graded algebra $({\cal M},d)$, where ${\cal M}$ is the free algebra
${\cal M}=\bigwedge (a_1, a_2, a_3,a_4)$ with all the generators of
degree $1$, and $d$ is given by $da_i=0$ for $i=1,2,4$ and
$da_3=-a_1 \cdot a_2$. The morphism $\rho\colon {\cal M} \to \Omega(KT)$,
inducing an isomorphism on cohomology, is defined by
$\rho(a_1)=\alpha$, $\rho(a_2)=\beta$, $\rho(a_3)=\gamma$,
$\rho(a_4)=\eta$.

Now, according to Definition~\ref{primera}, $C^1=\la
a_1,a_2,a_4\ra$ and $N^1=\la a_3\ra$. Since the element $a_1 \cdot
a_3 \in N^1\cdot V^1$ is closed but not exact, we conclude that
$({\cal M},d)$ is not $1$--formal, and by Theorem~\ref{criterio2}
it is not formal. This fact is also a consequence of
Lemma~\ref{formalitynilm}. (The non-formality of $KT$ was proved
in~\cite{topo} seeing the existence of non-trivial Massey
products, and in~\cite{Ha} proving that tori are the only compact
nilmanifolds with a formal minimal model.)

A symplectic form $\omega$ on $KT$ is $\omega=\beta \wedge \gamma
+ \alpha \wedge \eta$. Since $[\omega] \cup [\beta]=0$ in
$H^3(KT)$, we get that $KT$ does not have the $1$--Lefschetz
property. (In~\cite{BG} it is proved that tori are the only
compact symplectic nilmanifolds satisfying the $1$--Lefschetz
property.)

The Kodaira-Thurston manifold can be also defined as a
$T^{2}$--bundle over $T^{2}$~\cite{Th}, and the symplectic form
$\omega$ defines an integral cohomology class. It is clear that
any Donaldson submanifold $Z$ of $KT$ is a symplectic manifold of
dimension $2$. Hence {\em Z is a K\"ahler manifold} and thus
formal and hard Lefschetz. This result is true in general for the
Donaldson submanifolds of any $4$--dimensional compact symplectic
manifold.

For simplicity we shall denote by the same symbols the
differential forms on $KT$ and the ones induced in $Z$. {}From
Corollary~\ref{cohomologyDonaldson} it follows that the cohomology
classes  $[\alpha\wedge \gamma]$ and $[\beta \wedge \eta]$ of $KT$
define the zero class in $Z$. Moreover, $[\alpha \wedge \eta]$ and
$[\beta \wedge \gamma]$ restrict to the same cohomology class in
$Z$. Therefore, the real cohomology of $Z$ is
\begin{eqnarray*}
 H^0(Z) &=& \la 1\ra, \\
 H^1(Z) &=& \la [\alpha], [\beta], [\eta], [e_k]\ra,\\
 H^2(Z) &=& \la [\alpha \wedge \eta]\ra,
\end{eqnarray*}
where $[e_k]$ are a finite number of cohomology classes lying in $Z$ but
not in $KT$.

{\bf Example 2} {\it The Iwasawa manifold\/}. Consider the complex
Heisenberg group, that is, the complex nilpotent Lie group $G$ of
complex matrices of the form
 $$
 \pmatrix{1&x&z\cr 0&1&y\cr 0&0&1 \cr}.
 $$
The Iwasawa manifold is the compact complex parallelizable
nilmanifold obtained as $I_3=\Gamma{\backslash}G$, where $\Gamma$
is the discrete subgroup of $G$ consisting of those matrices whose entries
are Gaussian integers. The (complex) differential forms $dx$,
$dy$, $dz-xdy$ on $G$ are left invariant and descend to
holomorphic $1$--forms $\alpha$, $\beta$, $\gamma$ on $I_3$ such
that
 $$
 d\alpha = d\beta =0,  \quad  d\gamma = -\alpha\wedge\beta.
 $$
Denote by $\alpha_1=\Re(\alpha)$, $\alpha_2=\Im(\alpha)$,
$\beta_1=\Re(\beta)$, $\beta_2=\Im(\beta)$,
$\gamma_1=\Re(\gamma)$, $\gamma_2=\Im(\gamma).$ Then, $\{\alpha_1,
\alpha_2, \beta_1, \beta_2, \gamma_1, \gamma_2\}$ is a basis for
the $1$--forms on $I_3$ such that
 \begin{eqnarray*}
 d\alpha_i &= &d \beta_i =0,\\
 d\gamma_1 &=& -\alpha_1 \wedge \beta_1 + \alpha_2 \wedge \beta_2,\\
 d\gamma_2 &=& -\alpha_1 \wedge \beta_2 - \alpha_2 \wedge \beta_1.
 \end{eqnarray*}
Nomizu's theorem~\cite{No} implies that the real cohomology of $I_3$ is
 \begin{eqnarray*}
 H^0(I_3) &=& \la 1\ra,\\
 H^1(I_3) &=& \la [\alpha_1], [\alpha_2], [\beta_1], [\beta_2]\ra,\\
 H^2(I_3) &=& \la [\alpha_1 \wedge \alpha_2], [\alpha_1 \wedge \beta_1],
  [\alpha_1 \wedge \beta_2], [\beta_1\wedge \beta_2],
 [\alpha_1\wedge \gamma_2 + \alpha_2\wedge \gamma_1],
 [\alpha_1\wedge \gamma_1 - \alpha_2\wedge \gamma_2],\\
  & &
  [\beta_1\wedge \gamma_2 + \beta_2\wedge \gamma_1],
  [\beta_1\wedge \gamma_1 - \beta_2\wedge \gamma_2]\ra,\\
 H^3(I_3) &=& \la [\alpha_1 \wedge \alpha_2 \wedge \gamma_1],
  [\alpha_1 \wedge \alpha_2 \wedge \gamma_2],
  [\beta_1 \wedge \beta_2 \wedge \gamma_1],
  [\beta_1 \wedge \beta_2 \wedge \gamma_2],
  [\alpha_1 \wedge \beta_1 \wedge \gamma_2],\\
 & &
   [\alpha_1 \wedge \beta_2 \wedge \gamma_1],
  [\alpha_2 \wedge \beta_1 \wedge \gamma_1],
  [\alpha_1 \wedge \beta_1 \wedge \gamma_1-\alpha_1 \wedge \beta_2 \wedge
     \gamma_2], \\
 & &
  [\alpha_1 \wedge \beta_2 \wedge \gamma_2-\alpha_2 \wedge \beta_1 \wedge
     \gamma_2],
  [\alpha_2 \wedge \beta_1 \wedge \gamma_2+\alpha_2 \wedge \beta_2 \wedge
     \gamma_1]\ra,\\
 H^4(I_3) &=& \la [\alpha_1 \wedge \alpha_2 \wedge \beta_1 \wedge \gamma_1],
   [\alpha_1 \wedge \alpha_2 \wedge \beta_1 \wedge \gamma_2],
   [\alpha_1 \wedge \alpha_2 \wedge \gamma_1 \wedge \gamma_2],
   [\alpha_1 \wedge \beta_1 \wedge \beta_2 \wedge \gamma_1],\\
  & &   [\alpha_1 \wedge \beta_1 \wedge \beta_2 \wedge \gamma_2],
   [\beta_1 \wedge \beta_2 \wedge \gamma_1 \wedge \gamma_2],
   [\alpha_1 \wedge \beta_1 \wedge \gamma_1 \wedge \gamma_2+
    \alpha_2 \wedge \beta_2 \wedge \gamma_1 \wedge \gamma_2],\\
  & & [\alpha_1 \wedge \beta_2 \wedge \gamma_1 \wedge \gamma_2-
    \alpha_2 \wedge \beta_1 \wedge \gamma_1 \wedge \gamma_2]\ra,\\
 H^5(I_3) &=& \la [\alpha_1 \wedge \alpha_2 \wedge \beta_1 \wedge
  \gamma_1 \wedge \gamma_2],
  [\alpha_1 \wedge \alpha_2 \wedge \beta_2 \wedge
  \gamma_1 \wedge \gamma_2],
  [\alpha_1 \wedge \beta_1 \wedge \beta_2 \wedge
  \gamma_1 \wedge \gamma_2],\\
 & &
  [\alpha_2 \wedge \beta_1 \wedge \beta_2 \wedge
  \gamma_1 \wedge \gamma_2]\ra, \\
 H^6(I_3) &=& \la [\alpha_1 \wedge \alpha_2 \wedge \beta_1 \wedge \beta_2
  \wedge \gamma_1 \wedge \gamma_2]\ra.
 \end{eqnarray*}
{}From Lemma~\ref{formalitynilm} we know that $I_3$ is not
$1$--formal. Independently, one can check that the minimal model
of $I_3$ is the differential graded algebra $({\cal M},d)$, where
${\cal M}$ is the free algebra ${\cal M}=\bigwedge (a_1, a_2, b_1,
b_2, c_1, c_2)$ with all the generators of degree $1$, and $d$ is
given by $da_i=db_i=0$ for $1\leq i\leq 2$, $dc_1=-a_1 \cdot b_1 +
a_2 \cdot b_2$ and $dc_2=-a_1 \cdot b_2 - a_2 \cdot b_1$. The
morphism $\rho\colon {\cal M} \to \Omega(I_3)$, inducing an
isomorphism on cohomology, is defined by $\rho(a_i)=\alpha_i$,
$\rho(b_i)=\beta_i$, $\rho(c_i)=\gamma_i$, ($i=1, 2$).

Now, according to Definition~\ref{primera}, $C^1=\la
a_1,a_2,b_1,b_2\ra$ and $N^1=\la c_1,c_2\ra$. Since the element
$c_1 \cdot a_1 \cdot a_2$ in the ideal generated by $N^1$ in
$\bigwedge V^1$ is closed but not exact, we conclude that $({\cal
M},d)$ is not $1$--formal, and by Theorem~\ref{criterio2} it is
not formal. Therefore, $I_3$ is not $1$--formal, and thus
non-formal.

A symplectic form $\omega$ on $I_3$ is given by $\omega= \alpha_1
\wedge \gamma_2 + \alpha_2 \wedge \gamma_1 + \beta_1 \wedge
\beta_2$. It is easy to show that $[\omega]^2 \cup [\alpha_1]=0$,
so $I_3$ does not have the $1$--Lefschetz property.

The Iwasawa manifold can be also defined as a $T^{2}$--bundle over
$T^{4}$, and the symplectic form $\omega$ defines an integral
cohomology class. Let $Z \inc I_3$ be a Donaldson submanifold of
$I_3$. Then $Z$ is a $4$--dimensional symplectic manifold. For
simplicity we shall denote by the same symbols the differential
forms on $I_3$ and the ones induced in $Z$. Using
Corollary~\ref{cohomologyDonaldson} one can check that the real
cohomology of $Z$ is
\begin{eqnarray*}
 H^0(Z) &=& \la 1\ra, \\
 H^1(Z) &=& \la [\alpha_1], [\alpha_2], [\beta_1], [\beta_2]\ra,\\
 H^2(Z) &=& \la [\alpha_1 \wedge \alpha_2], [\alpha_1 \wedge \beta_1],
  [\alpha_1 \wedge \beta_2],
 [\beta_1\wedge \beta_2],
 [\alpha_1\wedge \gamma_2 + \alpha_2\wedge \gamma_1],
  [\alpha_1\wedge \gamma_1 - \alpha_2\wedge \gamma_2],\\
 & &
 [\beta_1\wedge \gamma_2 + \beta_2\wedge \gamma_1],
  [\beta_1\wedge \gamma_1 - \beta_2\wedge \gamma_2], [e_k]\ra,\\
 H^3(Z) &=& \la  [\beta_1 \wedge \beta_2 \wedge \gamma_1],
  [\beta_1 \wedge \beta_2 \wedge \gamma_2],
  [\alpha_1 \wedge \beta_1 \wedge \gamma_2],
  [\alpha_1 \wedge \beta_1 \wedge \gamma_1-\alpha_1 \wedge \beta_2 \wedge
     \gamma_2]\ra,\\
 H^4(Z) &=& \la [\alpha_1 \wedge \alpha_2 \wedge \gamma_1
   \wedge \gamma_2]\ra,
\end{eqnarray*}
where $[e_k]$ are a finite number of cohomology classes of $Z$ that are
not defined in $I_3$.

\begin{proposition} \label{Donaldson-Iwasawa}
Any Donaldson submanifold $Z \inc I_3$ is a $4$--dimensional
symplectic manifold not formal and not hard Lefschetz. Moreover,
$Z$ does not carry complex structures.
\end{proposition}

\begin{proof}
First let us note that, as a consequence of
Theorem~\ref{Lefschetz-Donaldson}, any Donaldson submanifold $Z$
of $I_3$ does not satisfy the hard Lefschetz theorem. By
Proposition~\ref{equivalence1}, the minimal model of $Z$ is of the
form $({\cal M}_Z,d)$, where ${\cal M}_Z=\bigwedge \left(
(a_1,a_2,b_1,b_2,c_1,c_2 ) \oplus V^{\geq 2}\right)$. Setting
$C^1=\la a_1,a_2,b_1,b_2\ra$ and $N^1=\la c_1,c_2\ra$, the element
$c_1 \cdot a_2 + c_2 \cdot a_1 \in N^1 \cdot \bigwedge V^1$ is
closed but not exact in $H^*(Z)$. Therefore by $Z$ is not
$1$--formal and hence it is not formal.

To show that $Z$ has no complex structures, we use Kodaira's
theorem~\cite{Ko} that states that {\it a complex surface is a
deformation of an algebraic surface if and only if its first Betti
number is even}. Suppose $Z$ with first Betti number $b_1(Z)=4$
has a complex structure. Then Kodaira's theorem implies that $Z$
possesses a K\"ahler metric, and hence $Z$ would be formal
according to a result of~\cite{DGMS}. But this is impossible.
\end{proof}

{\bf Example 3} {\it The manifold of~\cite{FLS}\/}. Let $G$ be the
connected completely solvable Lie group of dimension $6$
consisting of matrices of the form
 $$
 a=\pmatrix{e^t&0&xe^t&0&0&y_1 \cr
 0&e^{-t}&0&xe^{-t}&0&y_2\cr 0&0&e^t&0&0&z_1\cr 0&0&0&e^{-t}&0&z_2
 \cr 0&0&0&0&1&t \cr 0&0&0&0&0&1 \cr},
 $$
where $t, x, y_i, z_i \in \RR$ ($i=1,2$). Then a global system of
coordinates ${t,x,y_1,y_2,z_1,z_2}$ for $G$ is defined by
$t(a)=t$, $x(a)=x$, $y_i(a)=y_i$, $z_i(a)=z_i$, and a standard
calculation shows that a basis for the left invariant $1$--forms
on $G$ consists of
    $$
    \{dt,dx, e^{-t}dy_1-xe^{-t}dz_1, e^tdy_2-xe^tdz_1, e^{-t}dz_1, e^{-t}dz_2\}.
    $$
Let $\Gamma$ be a discrete subgroup of $G$ such that the quotient
space $M =\Gamma\backslash G$ is compact. (Such a subgroup exists,
see~\cite{FLS}.) Hence the forms $dt,dx, e^{-t}dy_1-xe^{-t}dz_1,
e^tdy_2-xe^tdz_1, e^{-t}dz_1, e^{-t}dz_2$  descend to $1$--forms
$\alpha$, $\beta$, $\gamma_1$, $\gamma_2$, $\delta_1$, $\delta_2$
on $M$ such that
 $$
 d\alpha=d\beta=0, \quad
 d\gamma_1=-\alpha \wedge \gamma_1 - \beta \wedge \delta_1,\quad
 d\gamma_2=\alpha \wedge \gamma_2 - \beta \wedge \delta_2,\quad
 d\delta_1=- \alpha \wedge \delta_1,\quad
 d\delta_2= \alpha \wedge \delta_2,
 $$
and such that at each point of $M$, $\{\alpha, \beta, \gamma_i,
\delta_i\}$ is a basis for the $1$--forms on $M$. Using Hattori's
theorem~\cite{Hat} we compute the real cohomology of $M$:
 \begin{eqnarray*}
 H^0(M) &=& \la 1\ra,\\
 H^1(M) &=& \la [\alpha], [\beta]\ra,\\
 H^2(M) &=& \la [\alpha \wedge \beta], [\delta_1\wedge \delta_2],
 [\gamma_1 \wedge \delta_2 + \gamma_2 \wedge \delta_1]\ra,\\
 H^3(M) &=& \la [\alpha \wedge \delta_1 \wedge \delta_2],
  [\beta\wedge \gamma_1 \wedge \gamma_2],
  [\beta\wedge (\gamma_1 \wedge \delta_2 + \gamma_2 \wedge \delta_1)],
  [\alpha\wedge (\gamma_1 \wedge \delta_2 + \gamma_2 \wedge
  \delta_1)]\ra,\\
 H^4(M) &=& \la [\alpha \wedge \beta \wedge \gamma_1 \wedge \gamma_2],
   [\alpha \wedge \beta \wedge \gamma_1 \wedge \delta_2],
   [\gamma_1 \wedge \gamma_2 \wedge \delta_1\wedge \delta_2]\ra,\\
 H^5(M) &=& \la [\alpha \wedge \gamma_1 \wedge \gamma_2 \wedge \delta_1\wedge
   \delta_2], [\beta \wedge \gamma_1 \wedge \gamma_2 \wedge \delta_1\wedge
   \delta_2]\ra, \\
 H^6(M) &=& \la [\alpha \wedge \beta \wedge \gamma_1 \wedge \gamma_2 \wedge
       \delta_1\wedge  \delta_2]\ra.
 \end{eqnarray*}
The minimal model of $M$ must be a differential graded algebra
$({\cal M},d)$, being ${\cal M}$ the free algebra of the form
${\cal M}=\bigwedge(a_1, a_2) \otimes \bigwedge(b_1,b_2,b_3,b_4)
\otimes \bigwedge V^{\geq 3}$ where the generators $a_i$ have
degree $1$, the generators $b_j$ have degree 2, and $d$ is given
by $da_i=db_1=db_2=0$, $db_3=-a_2 \cdot b_1$ and $db_4=a_2 \cdot
b_3$. The morphism $\rho\colon {\cal M} \to \Omega(M)$, inducing
an isomorphism on cohomology, is defined by $\rho(a_1)=\alpha$,
$\rho(a_2)=\beta$, $\rho(b_1)=\delta_1 \wedge \delta_2$,
$\rho(b_2)=(1/2)(\gamma_1 \wedge \delta_2 +\gamma_2 \wedge
\delta_1)$, $\rho(b_3)=(1/2)(\gamma_2 \wedge \delta_1 -\gamma_1
\wedge \delta_2)$ and $\rho(b_4)=(1/2) (\gamma_1 \wedge
\gamma_2)$.

According to Definition~\ref{primera}, we get $C^1=\la a_1,a_2\ra$
and $N^1=0$, thus {\em M is 1--formal}. We see that {\em $M$ is
not 2--formal}. In fact, the element $b_4 \cdot a_2 \in N^2 \cdot
V^1$ is closed but not exact, which implies that $({\cal M},d)$ is
not $2$--formal, and by Theorem~\ref{criterio2} not formal.

Consider the symplectic form $\omega$ on $M$ given by
$\omega=\alpha \wedge \beta + \gamma_1 \wedge \delta_2 + \gamma_2
\wedge \delta_1$. Then ($M,\omega$) is a {\em $1$--formal
symplectically aspherical manifold but not $2$--formal}.

Moreover, $[\omega] \cup [\delta_1 \wedge \delta_2]=0$ in
$H^4(M)$, which means that $M$ does not have the $2$--Lefschetz
property. But a simple computation shows that the cup product by
$[\omega]^2$ is an isomorphism between $H^1(M)$ and $H^5(M)$.
Therefore, $M$ has the $1$--Lefschetz property.


Alternatively, the manifold $M$ can be also defined as a
$T^{4}$--bundle over $T^{2}$ (see~\cite{FLS}), and the symplectic
form $\omega$ defines an integral cohomology class. Let $Z \inc
I_3$ be a Donaldson submanifold. For simplicity we shall denote by
the same symbols the differential forms on $M$ and the ones
induced in $Z$. {}From Corollary~\ref{cohomologyDonaldson} it
follows that the cohomology classes $[\alpha \wedge \delta_1
\wedge \delta_2]$, $[\beta \wedge \gamma_1 \wedge \gamma_2]$,
$[\alpha \wedge \beta \wedge \gamma_1 \wedge \gamma_2]$ and
$[\alpha \wedge \beta \wedge \gamma_1 \wedge \delta_2]- [\gamma_1
\wedge \gamma_2 \wedge \delta_1\wedge \delta_2]$ of $M$ define the
zero class in $Z$.
 Therefore, the real cohomology of $Z$ is
 \begin{eqnarray*}
 H^0(Z) &=& \la 1\ra,\\
 H^1(Z) &=& \la [\alpha], [\beta]\ra,\\
 H^2(Z) &=& \la [\alpha \wedge \beta], [\delta_1\wedge \delta_2],
 [\gamma_1 \wedge \delta_2 + \gamma_2 \wedge \delta_1], [e_k]\ra,\\
 H^3(Z) &=& \la [\beta\wedge (\gamma_1 \wedge \delta_2 + \gamma_2 \wedge
\delta_1)],
  [\alpha\wedge (\gamma_1 \wedge \delta_2 + \gamma_2 \wedge
  \delta_1)]\ra,\\
 H^4(Z) &=& \la [\alpha \wedge \beta \wedge \gamma_1 \wedge \delta_2]\ra,
 \end{eqnarray*}
where $[e_k]$ are a finite number of cohomology classes of $Z$ that
are not defined in $M$.

\begin{proposition}\label{Jaume}
Any Donaldson submanifold $Z \inc M$ is a $4$--dimensional formal
symplectic manifold that satisfies the hard Lefschetz property. Moreover,
$Z$ does not admit complex structures and, in particular, $Z$
does not possess K\"{a}hler metrics.
\end{proposition}

\begin{proof}
Since $M$ is $1$--formal and has the $1$--Lefschetz property, $Z$
is formal and hard Lefschetz. Suppose that $Z$ has no K\"{a}hler
metrics. Using Kodaira's theorem and $b_1(Z)=2$, a similar
argument to the one given in Proposition~\ref{Donaldson-Iwasawa}
implies that $Z$ has no complex structures if $Z$ has no K\"ahler
metrics.

In order to show that $Z$ does not admit K\"{a}hler metrics, recall
that $\Gamma=\pi_1(Z) \iso \pi_1(M)$, which is a semi-direct
product $\ZZ^2 \ltimes \ZZ^4$, so $\Gamma$ is $2$--step solvable.
Moreover, its rank is is $6$ by additivity. We shall see that
$\Gamma$ cannot be the fundamental group of any compact K\"{a}hler
manifold.

Assume now that $\Gamma= \pi_1(X)$, where $X$ is a compact K\"{a}hler
manifold. According to Arapura--Nori's theorem (see Theorem 3.3
of~\cite{AN}), there exists a chain of normal subgroups
 $$
 0 = D^3 \Gamma \subset Q \subset P \subset \Gamma,
 $$
such that $Q$ is torsion, $P/Q$ is nilpotent and $\Gamma/P$ is
finite. Since $\Gamma$ has no torsion, $Q=0$. As $\Gamma/ P$ is
torsion, we have $\rank P= \rank \Gamma= 6$. Now, the finite
inclusion $P \subset \Gamma$ defines a finite cover $p: Y
\rightarrow X$ that is also compact K\"{a}hler and it has fundamental
group $P$. By Corollary 3.8 of~\cite{Ca}, as $P$ is K\"ahler,
nilpotent and has $\rank P=6<9$, it has to be abelian. This is
impossible since any pair of non-zero elements $e\in \ZZ^2\subset
\Gamma=\ZZ^2 \ltimes \ZZ^4$, $f \in \ZZ^4\subset\Gamma$ do not
commute (see \cite[page 22]{FLS}).
\end{proof}

\begin{remark}
In~\cite{FLS} it is proved that $M$ is not formal showing that
the quadruple Massey product $\la [\delta_1\wedge  \delta_2],
[\beta], [\beta], [\beta]\ra$ is non-trivial since any
representative is of the form $[\beta \wedge \gamma_1 \wedge
\gamma_2] + [\delta_1\wedge \delta_2] \cup [u_1] + [\beta] \cup
[u_2] + \lambda [\beta\wedge (\gamma_1 \wedge \delta_2 + \gamma_2
\wedge \delta_1)]$, where $\lambda$ is a real number, and $[u_i]
\in H^i(M)$ for $i=1,2$. In particular, a representative of that
Massey product is $[\beta \wedge \gamma_1 \wedge \gamma_2]$. Now,
let $Z \inc M$ be a Donaldson submanifold of ${(M,\omega)}$. Using
Corollary~\ref{cohomologyDonaldson}, one can check that $[\beta
\wedge \gamma_1 \wedge \gamma_2]$  defines the zero class in
$H^3(Z)$, and so the quadruple Massey product $\la [\delta_1\wedge
\delta_2], [\beta], [\beta], [\beta]\ra$ vanishes in $Z$.
Moreover, the Massey product $\la [\alpha], [\alpha], [\delta_1
\wedge  \delta_2]\ra$ is trivial in $Z$ but it is not defined in
$M$.
\end{remark}
\medskip

{\bf Example 4} {\it The manifold $V$\/}. Let us show an example
of a Donaldson symplectic submanifold which is formal and not hard
Lefschetz.

In~\cite{IRTU} the authors present an example of a simply
connected compact symplectic manifold $V$ with all Massey products
(of all orders) vanishing, and such that $V$ does not satisfy the
hard Lefschetz theorem. Now, our purpose is to prove that the
manifold $V$ is formal.

Recall shortly the definition of $V$. Let $(M,\omega)$ be a
$4$--dimensional compact symplectic manifold whose first Betti
number is $b_1(M)=1$. Such a manifold exists as consequence of the
results of Gompf~\cite{Go}. Without loss of generality we can
assume  that the symplectic form on $M$ is integral and therefore,
by Gromov and Tischler theorem~\cite{Gr1, Ti}, there exists a
symplectic embedding of $M$ in the complex projective space
$\CP^5$ with its standard K\"ahler form.

Denote by $X$ the blow up of $\CP^5$ along $M$. Then $X$ is a
simply connected compact symplectic manifold of dimension $10$
whose third Betti number $b_3(X)=1$. Define $V$ as a Donaldson
symplectic submanifold of $X$. Then $V$ is an eight dimensional
simply connected compact symplectic manifold. For $i<4$ the de
Rham cohomology groups $H^i(X)$ and $H^i(V)$ are isomorphic, and
there is a monomorphism $H^4(X)\inc H^4(V)$.

\begin{lemma}~\cite{IRTU}
The de Rham cohomology group $H^2(X)$ is generated by two elements
$\rho$ and $\sigma$ satisfying that the cup product $\rho^2 \cup
\sigma^2$ is a nonzero cohomology class in $H^8(X)$.
\end{lemma}

\begin{theorem} \label{formal-notLefschetz}
$X$ is $3$--formal but not $4$--formal, and it is not
$3$--Lefschetz. Therefore $V$ is formal and not hard Lefschetz.
\end{theorem}

\begin{proof}
Since $H^3(X)$ and $H^3(V)$ are isomorphic, $b_3(X)=b_3(V)=1$. It
is well known and easy to see that $b_{2i+1}(M)$ is even for any
hard Lefschetz $2n$--manifold $M$, since otherwise the product
bilinear form $H^{2i+1}(M)\ox H^{2n-2i-1}(M) \to H^{2n}(M)$ is
degenerate (see for example~\cite{IRTU}). So, $X$ and $V$ are not
hard Lefschetz.

From Corollary~\ref{8-manifolds}, we know that $X$ and $V$ are
$2$--formal. To prove that $V$ is formal, first we show that $X$
is $3$--formal. It is clear that the three cohomology classes
$\rho^2$, $\sigma^2$ and $\rho \cup \sigma$ must be non-trivial
because $\rho^2 \cup \sigma^2$ is a non-trivial class. This means
that, at the level of the minimal model $(\bigwedge V_X,d)$ of
$X$, the subspace $N^3$ (consisting of the non-closed generators
of $V_X$ of degree 3) is the zero space. Then, taking into account
Definition~\ref{primera}, we get that $(\bigwedge V_X,d)$ is
$3$--formal and so $X$ is $3$--formal. Now, we obtain the
formality of $V$ from Theorem~\ref{equivalence2}.
\end{proof}

\medskip

{\bf Example 5} {\it Manifolds which are $s$--formal but not
$(s+1)$--formal for any $s$\/}. To finish this section, we are
going to show that the notion of $s$--formality is not vacuous, by
giving some examples of compact oriented manifolds which are
$s$--formal but not $(s+1)$--formal, for any value of $s\geq 0$.
First note that Example 1 covers the case $s=0$ and Example 3
covers the case $s=1$ (in the case $s=1$, the manifold must be
non-simply connected, since otherwise Theorem \ref{Miller} implies
that $M$ is $2$--formal). So we restrict to $s\geq 2$. (By the
way, one can check the manifold $X$ defined in Example 4 is
$3$--formal but not $4$--formal.)

The examples that we are going to construct follow the pattern of
\cite{Ru}. They are not symplectic; but maybe with a little bit
more work one could obtain symplectic examples.

First we deal with the case $s=2i$ even. Consider a wedge of three
spheres $S^{i+1}\vee S^{i+1}\vee S^{i+1} \subset {\RR}^{i+2}$ and
let $\gamma \in \pi_{3i+1}(S^{i+1}\vee S^{i+1}\vee S^{i+1})$ be
the element represented by the iterated Whitehead product $\left[
[\iota_1,\iota_2],\iota_3\right]$ where $\iota_j$ is the image of
the generator of $\pi_{i+1}(S^{i+1})$ for the inclusion of
$S^{i+1}$ as the $j$--th factor in $S^{i+1}\vee S^{i+1}\vee
S^{i+1}$. Let $C$ be the cone of $\gamma: S^{3i+1} \too
S^{i+1}\vee S^{i+1}\vee S^{i+1}$. By \cite{Ru} there is a PL
embedding $C\subset {\RR}^{(3i+1)+(i+2)+1} = {\RR}^{4i+3}$. Note
that $C$ is $i$--connected and that the only non-vanishing
homology groups are $H^{i+1}(C)$ and $H^{3i+2}(C)$. Let $\alpha_j
\in H^2(C)$ be the cohomology class evaluating $1$ on the $j$--th
copy of $S^{i+1}$ and $0$ on the other two. Therefore the Massey
product $\la \a_1,\a_2,\a_3\ra \in H^{3i+2}(C)/([\a_1]\cup
H^{2i+1}(C) + H^{2i+1}(C) \cup [\a_3]) = H^{3i+2}(C)$ is non-zero
(see Lemma 7 in \cite{UM}).

Let $W$ be a closed regular neighborhood of $C$ and let $Z= \bd W$
be its boundary. We can arrange easily that $Z$ is a {\em
smooth\/} manifold of dimension $4i+2$. There is an exact sequence
  $$
   \cdots \to H^{k}(C)=H^k(W) \to H^k(Z) \to H^{k+1}(W,Z)=H_{4i+3-(k+1)}(W)
   = H_{4i+2-k}(C) \to \cdots
  $$
This implies that $H^k(Z)=0$ for $k\leq i$. Also
$\pi_1(Z)=\pi_1(W-C)\iso \pi_1(W) =\pi_1(C)=1$. Therefore $Z$ is
$i$--connected. By Theorem \ref{Miller}, $Z$ is $2i$--formal. Let
us prove that $Z$ is not $(2i+1)$--formal. If it were then by
Theorem \ref{criterio2} it would be formal. Let us see that $Z$
has a non-vanishing Massey product. Let $j:Z\inc W$ be the
inclusion. Now use that $j_*:H^{3i+2}(C)=H^{3i+2}(W) \to
H^{3i+2}(Z)$ is injective (because $H_{i}(C)=0$) and
$H^{2i+1}(Z)=0$. It follows that the Massey product $\la
j^*\a_1,j^*\a_2,j^*\a_3\ra \in  H^{3i+2}(Z)$ is non-zero.

\medskip

To cover the case $s=2i-1$, $i\geq 2$, we start with $S^{i}\vee
S^{i+1}\vee S^{i+1} \subset {\RR}^{i+2}$ let $\gamma = \left[
[\iota_1,\iota_2],\iota_3\right]\in \pi_{3i}(S^{i}\vee S^{i+1}\vee
S^{i+1})$, where $\iota_1$ is the image of $\pi_{i}(S^{i})$ and
$\iota_2$, $\iota_3$ are the images of the generators of
$\pi_{i+1}(S^{i+1})$ for the second and third factors. Let
$C\subset  {\RR}^{4i+2}$ be the cone of $\gamma$ as before, and
$\alpha_1 \in H^i(C)$, $\alpha_2,\alpha_3 \in H^{i+1}(C)$ be the
obvious cohomology classes dual to $\iota_j$. The Massey product
$\la \a_1,\a_2,\a_3\ra \in H^{3i+1}(C)/([\a_1]\cup H^{2i+1}(C) +
H^{2i}(C) \cup [\a_3]) = H^{3i+1}(C)$ is non-zero.

Again, let $W$ be a closed regular neighborhood of $C$ and let $Z=
\bd W$ be its boundary, which is taken to be a smooth manifold of
dimension $4i+1$. Let $j:Z\inc W$ be the inclusion. Then the
Massey product $\la j^*\a_1,j^*\a_2,j^*\a_3\ra \in  H^{3i+1}(Z)$
is non-zero, as before. So $Z$ is not $2i$--formal. On the other
hand, $Z$ is $(i-1)$-connected and $H^i(Z)={\RR}$. So the minimal
model ${\mathcal{M}}_Z$ of $Z$ has $V^{<i}=0$ and $V^i={\RR}$. As
in the proof of Theorem \ref{Miller} this implies that
$N^{<(2i-1)}=0$. Let $\xi$ be a generator of $V^i$. If $i$ is odd
then $\xi\cdot\xi=0$ so also $N^{2i-1}=0$. If $i$ is even and
$\xi\cdot \xi$ is not exact then again $N^{2i-1}=0$. If $i$ is
even and $[\xi\cdot\xi]=0$ then $N^{2i-1}$ is $1$--dimensional and
generated by an element $v$. Any element in $I_{2i-1}(N^{\leq
(2i-1)})$ is of the form $v\cdot z$, with $z \in \bV^{\leq
(2i-1)}$. This cannot be closed unless $z=0$. So $Z$ is
$(2i-1)$--formal.

\bigskip

\noindent {\bf Acknowledgments.} We thank Jaume Amor\'os for
providing us with a proof of Proposition~\ref{Jaume}. We are
grateful to Greg Lupton for pointing us to \cite{St} and to Ignasi
Mundet, Aniceto Murillo and Daniel Tanr\'e for useful
conversations. We also thank Haynes Miller for useful advice
which improved the presentation of this work as well as the
anonymous referee for clarifying remarks,
and for asking us to include examples of manifolds which are
$s$--formal and not $(s+1)$--formal. This work has been partially
supported through grants DGICYT (Spain) Project
PB-97-0504-C02-01/02, BFM2001-3778-C03-02 and BFM2000-0024. Also
partially supported by The European Contract Human Potential
Programme, Research Training Network HPRN-CT-2000-00101.


{\small

\vspace{0.15cm}

\noindent{\sf M. Fern\'andez:} Departamento de Matem\'aticas,
Facultad de Ciencias, Universidad del Pa\'{\i}s Vasco, Apartado 644,
48080 Bilbao, Spain. {\sl E-mail:} mtpferol@lg.ehu.es

\vspace{0.15cm}

\noindent{\sf V. Mu\~noz:} Departamento de Matem\'aticas, Facultad
de Ciencias, Universidad Aut\'onoma de Madrid, 28049 Madrid,
Spain. {\sl E-mail:} vicente.munoz@uam.es}

\end{document}